\documentclass[a4paper,article,11pt,fleqn,oneside]{memoir} 
\setsecnumdepth{section}
\usepackage {amsmath,amsthm,amsxtra,amssymb}
\usepackage[fleqn,tbtags]{mathtools}      
\usepackage{upgreek} 

\usepackage[utf8]{inputenc}
\usepackage[english]{babel}
\usepackage{lmodern}   
\usepackage[T1]{fontenc}	 
\usepackage {bm} 
\usepackage{fix-cm} 
\usepackage{microtype}	
\usepackage{float} 
\usepackage[compatibility=false]{caption}

\checkandfixthelayout[nearest]

\usepackage[babel=true]{csquotes}
\usepackage{graphicx}
\usepackage[ansinew]{} 
\usepackage{fixltx2e,mparhack}
\usepackage[fixpdftex]{xcolor} 
\usepackage{tikz}
\usepackage{pgf}
\usetikzlibrary{fit} 
\usetikzlibrary[shapes.geometric]
\usetikzlibrary{matrix,arrows,petri,topaths,positioning,automata}

\setlrmarginsandblock{3.14cm}{3.5cm}{*}
\setulmarginsandblock{3.14cm}{3.3cm}{*}
\checkandfixthelayout

\headstyles{memman}
\chapterstyle{article}
\setsubsecheadstyle{\centering\normalfont}
\setbeforesubsecskip{-1\onelineskip plus -1ex minus -.2ex}
\setaftersubsecskip{1\onelineskip plus .2ex}
\usepackage[colorlinks=true, pdfstartview=FitV, linkcolor=blue, citecolor=blue, urlcolor=blue]{hyperref}

\let \rho   = \varrho
\let \phi   = \varphi
\let \leq   = \leqslant
\let \Phi = \varPhi
\let \Psi = \varPsi
\let \Gamma = \varGamma
\let \Theta = \varTheta
\let \triangle= \vartriangle

\DeclareMathOperator{\Set}{\textbf{Set}}

\DeclareMathOperator{\ke}{ker}
\DeclareMathOperator{\Sub}{Sub}
\DeclareMathOperator{\id}{id}
\DeclareMathOperator{\FGR}{\textbf{Graph}_{\textbf{\F}}}
\DeclareMathOperator{\FGROne}{\textbf{Graph}_{\textbf{\F}_\textbf{1}}}
\DeclareMathOperator{\FGRTwo}{\textbf{Graph}_{\textbf{\F}_\textbf{2}}}

\DeclareMathOperator{\pb}{pb}

\DeclareMathOperator{\col}{Col}

\DeclareMathOperator{\Cc}{C}
\DeclareMathOperator{\Uu}{U}
\DeclareMathOperator{\C}{C}
\DeclareMathOperator{\Hom}{H}
\DeclareMathOperator{\HomM}{H^{-}}
\DeclareMathOperator{\Sim}{S}
\DeclareMathOperator{\Graph}{Graph}
\DeclareMathOperator{\Pat}{Pat}
\DeclareMathOperator{\Imp}{Imp}

\newcommand {\F}{F} 
\newcommand {\G}{G} 
\newcommand {\Kk}{\mathcal{K}} 

\newcommand {\Tt}{\mathcal{T}}

\newcommand {\evr}{(E,V,g)}
\newcommand {\evg}{(E,V,g)}
\newcommand {\Pp}{\mathcal{P}}
\DeclareMathOperator{\Ob}{Ob}
\DeclareMathOperator{\Mor}{Mor}

\newcommand{\ug}[1][]{#1_{\scriptscriptstyle{g}}}
\newcommand{\vv}[1][]{#1_{\scriptscriptstyle{v}}}
\newcommand{\ee}[1][]{#1_{\scriptscriptstyle{e}}}
\newcommand{\xx}[1][]{#1_{\scriptscriptstyle{X}}}

\newcommand{\ui}[1][]{#1_{i\in I}}

\newcommand{\inkl}[2]{\iota_{(\scriptstyle{#1}{,}{\scriptstyle{#2}})}}
\newcommand{\forg}[2]{\left|{#1}\right| _{\scriptscriptstyle{#2}}}

\newcommand {\he}{^{\scriptscriptstyle{(1)}}}   
\newcommand {\hz}{^{\scriptscriptstyle{(2)}}}
\newcommand {\hd}{^{\scriptscriptstyle{(3)}}}
\newcommand {\hr}{^{\scriptscriptstyle{(R)}}}
\newcommand {\hP}{^{\scriptscriptstyle{(P)}}}
\newcommand {\hEps}{^{\scriptstyle{(\epsilon)}}}

\newcommand {\hi}{^{\scriptscriptstyle{(i)}}}
\newcommand {\hj}{^{\scriptscriptstyle{(j)}}}
\newcommand {\hk}{^{\scriptscriptstyle{(k)}}}
\newcommand {\hme}{^{\scriptscriptstyle{-}}}

\newcommand {\htheta}{^{(\theta)}}
\newcommand{\sub}[1][]{#1_{\scriptscriptstyle{0}}}

\newcommand{\mdw}[1]{\texttt{\textcolor{weis}{[#1]}}}
\newcommand{\ink}[1]{(#1)}
\definecolor{rot}{rgb}{1,0,0}
\definecolor{weis}{rgb}{1,1,1}

\theoremstyle{definition}
\newtheorem{defi}{Definition}[section]
\newtheorem{numbsp}[defi]{Example}
\theoremstyle{remark}
\newtheorem{numbem}[defi]{Remark}

\theoremstyle{plain}
\newtheorem{lemma}[defi]{Lemma}
\newtheorem{satz}[defi]{Theorem}
\newtheorem{prop}[defi]{Proposition}
\newtheorem{theorem}[defi]{Theorem}
\newtheorem{folg}[defi]{Corollary}
\newtheorem{hsatz} [defi] {Lemma}

\begin{document}

\title{A Unified Categorical Approach To Graphs}
\author{Christian J\"akel, \textit{christian.jaekel(at)tu-dresden.de }\\
Technische Universit\"at Dresden}
\date{July 2015}
\maketitle
\begin{abstract}
For a set-endofunctor $\F$, we extend the notion of universal $\F$-coalgebras to $\F$-graphs. These generalized coalgebras are models for various types of graphs, such as (un)directed (hyper)graphs, relational structures or fuzzy graphs. The induced morphisms coincide with graph homomorphisms. From this point of view, graphs are \enquote{co-like} structures and share features of universal coalgebras.\\ 
In this article, we explore the coalgebraic character of graphs and transfer coalgebraic concepts like cofreeness, simulations or Co-Birkhoff theorems to $\F$-graphs. Products and cofree constructions for $\F$-graphs turn out to be less restrictive than their coalgebraic counterparts. 
\end{abstract}

\section{Introduction}

The theory of $\F$-coalgebras is a well developed subject and arose from duality to universal algebra. For a set-endofunctor $\F:\Set\rightarrow\Set$, a coalgebra is a pair $(A,\alpha)$, consisting of a set $A$ and a structure map $\alpha: A\rightarrow\F A$. An introduction to coalgebras over $\Set$ can be found in \cite{GummRand} or \cite{Ruth}. The more general theory of coalgebras over arbitrary categories is treated in \cite{HughesDiss}. Furthermore, we assume the reader to be familiar with basic category theory. Almost all used concepts can be found in \cite{JoyCat} or \cite{borceux2008handbook}. Co-Birkhoff type theorems are treated in \cite{HughesCoBikhoff}.\\
Coalgebras provide a general framework for modeling transition systems, automata or topological spaces; although in the latter case, morphisms are not the usually desired continuous maps. The same problem arises in modeling directed graphs as transition structures, \textit{i.e.}, morphisms do not preserve source and target nodes, but rather neighborhoods of each node. One way of dealing with this issue is to consider weak homomorphisms, as for example done in \cite{WeakHomDiss}.\\
We chose another approach, by generalizing $F$-coalgebras to $F$-graphs over $\Set\times\Set$. This yields a triple $(E,V,g)$ with $g:E\rightarrow\F V$, such that the induced morphisms coincide with classical graph homomorphisms (see section \ref{GraphsHomSubStrSection}) and the functor $\F$ determines the graph's type. Thereby, we generalize previous works about the category of graphs, \textit{e.g.}, \cite{PH}, \cite{Williams} or \cite{DeCatGraph}.\\
A large part of the structural theory of $\F$-coalgebras transfers well to $\F$-graphs, leading to a unified model of graphs (factorizations and congruences: section \ref{FactCongDiagLSection}, colimits: section \ref{ColimLimSection}, characterization of iso-,epi- and monomorphisms: section \ref{IsoEpiMonoSection}, functors between categories of graphs: section  \ref{IsoEpiMonoSection} and conjunct sums: section \ref{CSumsGrTrafoSection}). Limits (section \ref{ColimLimSection}) and cofree objects (section \ref{CofteeFGraphs}), in contrast to their coalgebraic counterparts, exist without any condition on $\F$.\\
Graph homomorphisms, which are useful tools in graph theory (\cite{StrucAndSymm,NesHell}), are naturally related to categorical modeling. As one can expect, in our generalized setting, we will not present new results about specific types of graphs. Nevertheless, we introduce the diagram lemmas (section \ref{FactCongDiagLSection}) and Co-Birkhoff type theorems (section \ref{CoBirkSection}) to graphs. We also provide a categorical foundation for graph relations (section \ref{GraphrelSection}), which have been applied to undirected graphs in \cite{LongDiss}.\\


\section{
\texorpdfstring{$\F$}{F}-Graphs, \texorpdfstring{$\F$}{F}-Graph Homomorphisms And Substructures}\label{GraphsHomSubStrSection}
First, we will answer the question \enquote{What is a graph?} and define structure conform mappings. We will do this by extending the notion of coalgebras over $\Set$ (the category of sets and maps between sets) to $\F$-graphs over $\Set\times\Set$. Next, we will define and characterize substructures.
\subsection{\texorpdfstring{$\F$}{F}-Graphs And \texorpdfstring{$\F$}{F}-Graph Homomorphisms}

For the whole article, let $\F$ be a covariant endofunctor in $\Set$.
\begin{defi}[$\F$-Graph]
An \emph{$\F$-graph} or a \emph{graph of type $\F$} is a triple $\G=(E,V,g)$, consisting of an \emph{edge set} $E$, a \emph{vertex set} $V$, and a \emph{structure map} $g:E\rightarrow\F V$\footnote{The traditional order of $E$ and $V$ is interchanged, because the edge set is the domain of the structure map $g:E\rightarrow\F V$. Thus, the edge set $E$ is listed first.}. 
\end{defi}
We give some examples for possible choices of $\F$ and resultant graph structures. 
\begin{numbsp}
With the identity functor, we can model multisets by assigning an arity to each element of $E$. The functor $\F:V\mapsto \mathfrak{P}_{_{1,2}}V$, assigning to a set its singleton and two-element subsets, is a model for undirected graphs. Directed graphs can be modeled via $\F:V\mapsto V\times V$.\\ Let $\mathfrak{P}$ be the powerset functor. The functors $\F:V\mapsto \mathfrak{P}V$ or $\F:V\mapsto V\times \mathfrak{P}V$ represent hypergraphs or directed hypergraphs. For an arbitrary functor $\F$, we define colored or weighted graphs via $\tilde{\F}(V):=X\ee\times\F(X\vv\times V)$, where $X\ee$ represents a set of edge colors and $X\vv$, a set of node colors. For instance, fuzzy graphs are modeled with $X\ee=X\vv=[0,1]$. Additional structure can be added, \textit{e.g.}, through coloring with a monoid or lattice. Hybrid graphs can be described by taking the sum of different type functors, as for example, $\F V=(V\times V) + \mathfrak{P}_{_{1,2}} V $ for a model of graphs with directed and undirected edges. In the latter framework, relational systems fit in too. 
\end{numbsp}
From now on, if not stated otherwise, $\G$ will be an $\F$-graph and as such refer to a triple $\evr$.
\begin{defi}
If the structure map $g$ is injective, we call $\G$ \emph{simple}. For surjective maps $g$, the $\F$-graph $\G$ is called \emph{complete}.
\end{defi}

\begin{defi}[$\F$-Graph Homomorphism]\label{DefiHom}
Let two $\F$-graphs $\G\he=(E\he,V\he,g\he)$ and $\G\hz=(E\hz,V\hz,g\hz)$ be given. A \emph{homomorphism} $\phi$ from $\G\he $ to $\G\hz$ is a pair of maps $(\phi\ee,\phi\vv)$, where the \emph{edge map} $\phi\ee: E\he\rightarrow E\hz$ and the \emph{vertex map} $\phi\vv: V\he\rightarrow V\hz$ have to fulfill $g\hz\circ\phi\ee=\F\ink{\phi\vv}\circ g\he.$ This equation can conveniently be expressed in terms of the following commutative diagram.

\begin{figure}[H]
\begin{center}
\begin{tikzpicture}[description/.style={fill=white,inner sep=2pt},>=stealth']
\matrix (m) [matrix of math nodes, row sep=2.71em,
column sep=3.14em, text height=1.5ex, text depth=0.25ex]
{ E\he & E\hz \\
 \F V\he & \F V\hz. \\ };
\path[->,font=\scriptsize](m-1-1) edge node[auto] {$ \phi\ee $} (m-1-2);
\path[->,font=\scriptsize](m-1-1) edge node[auto,swap] {$ g\he $} (m-2-1);
\path[->,font=\scriptsize](m-2-1) edge node[auto,swap] {$ \F(\phi\vv) $} (m-2-2);
\path[->,font=\scriptsize](m-1-2) edge node[auto] {$ g\hz  $} (m-2-2);
\end{tikzpicture}
\end{center}
\end{figure}

A homomorphism is called \emph{injective (surjective, bijective)} if $\phi\ee$ and $\phi\vv$ are injective (surjective, bijective) maps.
\end{defi}

\begin{numbsp}
For $\F V=\mathfrak{P} V$, homomorphisms are incidence preserving maps. In case of $\F V=V\times V$, homomorphisms preserve source and target nodes. Considering colored graphs, we get color and structure preserving maps. 
\end{numbsp}

\begin{prop}
The class of all $\F$-graphs for a fixed type functor $\F$, together with $\F$-graph homomorphisms and componentwise composition, defines a category.
\begin{proof}
For any $\G$, we define the identity homomorphism $\id_{\G}:=(\id_{E},\id_{V}):\G\rightarrow\G$. Let $\phi:\G\he\rightarrow\G\hz$ and $\psi:\G\hz\rightarrow\G\hd$ be homomorphisms. Their composition $\psi\circ\phi$ is defined as $(\psi\ee\circ\phi\ee,\psi\vv\circ\phi\vv).$ It is straight forward to show that all axioms of a category are fulfilled. 
\end{proof}
\end{prop}
\begin{defi}
We will denote the category of $\F$-graphs and $\F$-graph homomorphisms by $\FGR$. The category is equipped with the obvious \emph{forgetful functor} $\Uu:\FGR\rightarrow\Set\times\Set$, which maps $\evr$ to $(E,V)$ and which is the identity on homomorphisms $\phi=(\phi\ee,\phi\vv)$. Hence, we consider $\FGR$ as a concrete category over $\Set\times\Set$ (see \cite{JoyCat}).
\end{defi}

\subsection{Substructures}

We define substructures by means of regular monomorphisms -- this are homomorphisms which occur as the equalizer of two parallel homomorphisms -- and thereby, as we will see in theorem \ref{RegularmonoCharacterization}, assure that substructures are always equipped with an injective embedding homomorphism. 
\begin{defi}\label{Teilgraph} Let $\G$ be an $\F$-graph. A \emph{subgraph} of $\G$ is an $\F$-graph $\G\sub=(E\sub,V\sub,g\sub)$, together with a regular monomorphism $\iota:=(\inkl{E\sub}{E},\inkl{V\sub}{V})\footnote{The map $\inkl{X}{Y}:X\hookrightarrow Y$ is the natural inclusion map for $X\subseteq Y$.}:\G\sub\rightarrowtail\G.$ Symbolically, we write $\G\sub\leq\G$.
\end{defi}

\begin{numbem}\label{bemeindeutig}
For an $\F$-graph $\G$, we choose $E\sub\subseteq E$ and $V\sub\subseteq V$. If a structure map $g\sub$ exists, such that the inclusion map 
\mbox{$\inkl{\G\sub}{\G}:\G\sub\hookrightarrow \G$} defines a homomorphism, it will follow from theorem \ref{RegularmonoCharacterization} that $\inkl{\G\sub}{\G}$ is a regular monomorphism.\\ In this case, the subgraph $\G\sub=(E\sub,V\sub,g\sub)$ is uniquely defined through $g\sub$. If there would be two structure maps $g\sub,\tilde{g}\sub:E\sub\rightarrow\F V\sub$, such that the canonical inclusion $\iota_{(\G\sub,\G)}:\G\sub\hookrightarrow \G$ defines a homomorphism, we would have: \[\F(\inkl{V\sub}{V})\circ g\sub=g\;\circ\inkl{E\sub}{E}=\F(\inkl{V\sub}{V})\circ\tilde{g}\sub.\]
Because $\inkl{V\sub}{V}$ is left cancellable, the homomorphism $\F(\inkl{V\sub}{V})$ is it too and $g\sub=\tilde{g}\sub$ follows. 
\end{numbem}

At this point, we cite two lemmas from the category $\Set$ (see \cite{GummRand}), which we will use now and later in this article.  

\begin{hsatz}[First Diagram Lemma For Sets]\label{DLFM}
Let $g:A\rightarrow C$ and $f:A\twoheadrightarrow B$ be mappings, where $f$ is surjective and $C\neq\emptyset$. There exists a unique map $h:B\rightarrow C$ with $h\circ f=g$ if and only if $\ke f\subseteq \ke g$\footnote{For $f:A\rightarrow B$, the \emph{kernel} of $f$ is defined as $\ke f:=\{(a,b)\in A\times A \mid a,b\in A\;\text{and}\;f(a)=f(b)\}.$}.
\end{hsatz}
\begin{hsatz}[Second Diagram Lemma For Sets]\label{SecDiagLemmaSet} 
Let $f:B\rightarrowtail A$ and $g: C\rightarrow A$ be maps, where $f$ is injective. There exists a unique map $h: C\rightarrow B$ with $f\circ h=g$ if and only if $g[C]\subseteq f[B]$\footnote{For $f:A\rightarrow B$, the image of a $A$ under $f$ is denoted as $f[A]$.}.
\end{hsatz}

With the second diagram lemma for sets (\ref{SecDiagLemmaSet}), we can characterize which subsets qualify to define a subgraph. The coalgebraic version of this lemma can be found in \cite[Lemma 4.4.]{GummRand}.
	
\begin{prop}\label{SubGraphChar} 
Let $\G=(E,V,g$) be given. We choose $E\sub\subseteq E$ and $V\sub\subseteq V$. A structure map $g\sub:E\sub\rightarrow\F V\sub$, such that $\G\sub=(E\sub,V\sub,g\sub)\leq\G$, exists if and only if for every $e\in E\sub$ there is some $\tilde{v}\in\F V\sub$ with $g(e)=\F(\inkl{V\sub}{V})(\tilde{v}).$
\begin{proof}
We apply lemma \ref{SecDiagLemmaSet}. A map $g\sub$ exists if and only if $(g\circ\;\inkl{E\sub}{E})[E\sub]\subseteq\F(\inkl{V\sub}{V})[\F V\sub]$.\qedhere 
\end{proof}
\end{prop}
\noindent
\begin{numbem}\label{standardfunctor}
A set endofunctor $\F:\Set\rightarrow\Set$ is called \emph{standard} if it preserves inclusions, \textit{i.e.}, $\F(\inkl{U}{V})=\inkl{\F U}{\F V}$ whenever $U\subseteq V$. If the type functor $\F$ is standard, the criterion from proposition \ref{SubGraphChar} simplifies to: $e\in E\sub \Rightarrow g(e)\in \F V\sub.$
\end{numbem}

\section{Factorizations, Congruences and Diagram Lemmas}\label{FactCongDiagLSection}

In this section, we will follow the presentation in \cite[Section 3.2]{GU}. First, we will show that every $\F$-graph homomorphism has a unique surjective-injective factorization. Next, we will treat congruence relations and generalize the diagram lemmas \ref{DLFM}, \ref{SecDiagLemmaSet} to $\FGR$. 

\subsection{Homomorphic Images And Factorizations}

Analogously to remark \ref{bemeindeutig}, the following can be shown:

\begin{prop}\label{surstruktur} Let $\phi:\G\he\twoheadrightarrow\G\hz$ be a surjective homomorphism. The $\F$-graph structure $g\hz$ on $\G\hz$ is uniquely determined by $\phi$ and $g\he$. \[g\hz:=\{(\phi\ee(e),\F(\phi\vv)(g\he(e)))\mid e\in E\he\}.\]
\end{prop}

\begin{defi}\label{hombild}
Let $\phi$ be a homomorphism from $\G\he$ to $\G\hz$. The \emph{image} of $\G\he$ with respect to $\phi$ is defined as $\phi[\G\he]:=(\phi\ee{[}E\he{]},\phi\vv{[}V\he{]})$ and the \emph{image restriction}, as $\phi':=(\phi\ee',\phi\vv'):\G\he\twoheadrightarrow\phi[\G\he]$. If $\phi$ is a surjective homomorphism, we call $\G\hz$ \emph{homomorphic image} of $\G\he$.
\end{defi}

\begin{lemma}\label{faktorisierung} Let $\phi:\G\he\rightarrow\G\hd$ be a homomorphism and $\phi=\chi\circ\psi$ a factorization with a surjective map $\psi:(E\he,V\he)\twoheadrightarrow(E\hz,V\hz)$ followed by an injective map $\chi:(E\hz,V\hz)\rightarrowtail(E\hd,V\hd)$. There exists a unique $\F$-graph structure $g\hz$ on $\G\hz$. With respect to $g\hz$, the maps $\phi $ and $\chi$ are homomorphisms.

\begin{center}
\begin{tikzpicture}[description/.style={fill=white,inner sep=2pt},>=stealth']
\matrix (m) [matrix of math nodes, row sep=3em,
column sep=2.5em, text height=1.5ex, text depth=0.25ex]
{ E\he & E\hz & E\hd\\
 \F V\he & \F V\hz & \F V\hd \\ };

\path[->>,font=\scriptsize](m-1-1) edge node[swap,auto] {$ \psi\ee$} (m-1-2);
\path[>->,font=\scriptsize](m-1-2) edge node[swap,auto] {$ \chi\ee $} (m-1-3);
\path[->,font=\scriptsize](m-1-3) edge node[auto] {$ g\hd $} (m-2-3);
\path[>->,font=\tiny](m-2-2) edge node[auto] {$ \F(\chi\vv) $} (m-2-3);
\path[->>,font=\tiny](m-2-1) edge node[auto] {$ \F(\psi\vv)  $} (m-2-2);
\path[->,font=\scriptsize](m-1-1) edge node[swap,auto] {$ g\he  $} (m-2-1);
\path[->,font=\scriptsize](m-1-1) edge [bend left=15]node[auto] {$ \phi\ee  $} (m-1-3);
\path[->,font=\scriptsize](m-2-1) edge [bend right=15] node[swap,auto] {$ \F(\phi\vv)  $} (m-2-3);
\path[dashed,->,font=\scriptsize](m-1-2) edge node[auto] {$ g\hz  $} (m-2-2);
\end{tikzpicture}
\end{center}

\begin{proof}
Because $\phi$ is a homomorphism, we have:
$g\hd\circ(\chi\ee\circ\psi\ee)=\F(\chi\vv\circ\psi\vv)\circ g\he=\F(\chi\vv)\circ\F(\psi\ee)\circ g\he.$
As $\chi\vv$ is injective, $\F(\chi\vv)$ is injective too. Hence, $\psi\ee$, $g\hd\circ\chi\ee$, $\F(\psi\vv)\circ g\he$ and $\F(\chi\vv)$ define an E-M-square in $\Set$. Thus, a unique diagonal exists, which is the desired structure map $g\hz$.
\end{proof}
\end{lemma}

\begin{satz}\label{hombildundteilgraph} Let $\phi:\G\he\rightarrow\G\hz$ be a homomorphism. Via $\phi=\inkl{\phi[\G\he]}{\G\hz}\circ\;\phi'$, we can define a surjective-injective factorization of $\phi$. Additionally, $\phi[\G\he]$ is a homomorphic image of $\G\he$ and a subgraph of $\G\hz$. 
\begin{proof}
Obviously, $\inkl{\phi[\G\he]}{\G\hz}$ is injective and $\phi'$ is surjective.

\begin{center}
\begin{tikzpicture}[description/.style={fill=white,inner sep=2pt},>=stealth']
\matrix (m) [matrix of math nodes, row sep=3.14em,
column sep=6.14em, text height=1.5ex, text depth=0.25ex]
{ E & \phi[E] & \tilde{E}\\
 \F V & \F\phi[V] & \F \tilde{V}\\ };

\path[->>,font=\scriptsize](m-1-1) edge node[swap,auto] {$ \phi'\ee$} (m-1-2);
\path[shorten >=2,>->,font=\scriptsize](m-1-2) edge node[swap,auto] {$ \inkl{\phi[E\he]}{E\hz} $} (m-1-3);
\path[->,font=\scriptsize](m-1-3) edge node[auto] {$ g\hz $} (m-2-3);
\path[shorten >=2,>->,font=\scriptsize](m-2-2) edge node[auto] {$ \F(\inkl{\phi[V\he]}{V\hz}) $} (m-2-3);
\path[->>,font=\scriptsize](m-2-1) edge node[auto] {$ \F(\phi'\vv)  $} (m-2-2);
\path[->,font=\scriptsize](m-1-1) edge node[swap,auto] {$ g\he  $} (m-2-1);
\path[shorten >=2,->,font=\scriptsize](m-1-1) edge [bend left=15]node[auto] {$ \phi\ee  $} (m-1-3);
\path[shorten >=2,->,font=\scriptsize](m-2-1) edge [bend right=15] node[swap,auto] {$ \F(\phi\vv)  $} (m-2-3);
\path[dashed,->,font=\scriptsize](m-1-2) edge node[auto] {$  $} (m-2-2);
\end{tikzpicture}
\end{center}

Because of lemma \ref{faktorisierung}, the maps $\phi'$ and $\inkl{\G\he}{\G\hz}$ are $\F$-graph homomorphisms. Hence, $\phi[\G\he]$ is the homomorphic image of $\G\he$, and $\phi[\G\he]$ a subgraph of $\G\hz$. \qedhere
\end{proof}
\end{satz}

If $\G\sub\he$ is a subgraph of $\G\he$ and $\phi:\G\he\rightarrow\G\hz$ a homomorphism, we can apply theorem \ref{hombildundteilgraph} to $\phi\,\circ\inkl{\G\sub\he}{\G\he}$ and get:

\begin{folg}
$\phi[\inkl{\G\sub\he}{\G\he}[\G\he\sub]]$ is a subgraph of $\G\hz$.
\end{folg}

\subsection{Congruence Relations And Diagram Lemmas}
We will define homomorphism kernels and congruences. Next, we will show the existence of factor graphs and proof the diagram lemmas for $\FGR$. 

\begin{defi} Let $\G\he$ and $\G\hz$ be $\F$-graphs and $\phi:\G\he\rightarrow\G\hz$ a homomorphism.
The \emph{kernel} of $\phi$ is defined as $\ke\phi:=(\ke\phi\ee,\ke\phi\vv)$. For any $\F$-graph, the kernels are partially ordered via: \[\ke\phi\subseteq\ke\psi:\Longleftrightarrow \ke\phi\ee\subseteq\ke\psi\ee\;\text{and}\;\ke\phi\vv\subseteq\ke\psi\vv.\]
A \emph{congruence} $\theta=(\theta\ee,\theta\vv)$ is a pair of equivalence relations $\theta\ee\subseteq E\times E$ and $\theta\vv\subseteq V\times V$, such that a homomorphism $\phi$ exists, with $\theta=\ke\phi$. If $\theta=(\theta\ee,\theta\vv)$ is a congruence on $\G$, then $\G/\theta:=(E/\theta\ee,V/\theta\vv,g\htheta)$ defines the \emph{factor graph} of $\theta$. The \emph{canonical map} $\pi_{\theta}:\G\twoheadrightarrow\G/\theta$ is defined as $\pi_{\theta}:=(\pi_{\theta\ee},\pi_{\theta\vv})$.
\end{defi}

The following theorem shows that a unique factor graph structure always exists.

\begin{satz}\label{Faktorgraph}
Let $\phi:\G\he\rightarrow\G\hz$ be a homomorphism and $\theta=(\theta\ee,\theta\vv)$ its kernel. For the factor graph $\G/\theta$, there exists a unique $\F$-graph structure $g\htheta$. 
\begin{proof}
The surjective canonical map $\pi_{\theta}:=(\pi_{\theta\ee},\pi_{\theta\vv})$ from $\G\he$ to $\G\he/\theta$ exists. Furthermore, we define $\psi:\G\he/\theta\rightarrow\phi[\G\he]$ with $\psi\ee:[e]_{\theta\ee}\mapsto\phi\ee(e)$ for $e\in E\he$ and $\psi\vv:[v]_{\theta\vv}\mapsto\phi\vv(v)$ for $v\in V\he$. It is easy to show that $\psi$ is well defined and injective. Consequently, there is a surjective-injective factorization $\phi'=\psi\circ\pi_{\theta}$. From lemma \ref{faktorisierung}, the existence and uniqueness of $g\htheta$ follows. 
\end{proof}
\end{satz}

\begin{defi} An equivalence relation on an $\F$-graph $\G=\evr$ is a pair $(R\ee,R\vv)$, where $R\ee$ and $R\vv$ are equivalence relations on $E$ and $V$ respectively. 
\end{defi}

Due to \cite[Corollar 3.2.8]{GU}, a criterion, whether a given equivalence relation is a congruence, is obtained by:

\begin{prop}\label{ConChar} An equivalence relation $\theta$ on $\G=(E,V,g)$ is a congruence if and only if	for all edges $e,\tilde{e}\in E$, we have that:
\[e\,\theta\ee \tilde{e}\Longrightarrow\F(\pi_{\theta\vv})(g(e))=\F(\pi_{\theta\vv})(g(\tilde{e})).\]
\begin{proof}
Every congruence $\theta$ is the kernel of $\pi_{\theta}$ and the following diagram commutes.

\begin{center}
\begin{tikzpicture}[description/.style={fill=white,inner sep=2pt},>=stealth']
\matrix (m) [matrix of math nodes, row sep=3em,
column sep=2.5em, text height=1.5ex, text depth=0.25ex]
{ E & E/{\theta\ee} \\
 \F V & \F V/\theta\vv \\ };
\path[->>,font=\scriptsize](m-1-1) edge node[auto] {$ \pi_{\theta\ee}$} (m-1-2);
\path[->,font=\scriptsize](m-1-1) edge node[auto,swap] {$ g $} (m-2-1);
\path[->>,font=\scriptsize](m-2-1) edge node[auto,swap] {$ \F(\pi_{\theta\vv}) $} (m-2-2);
\path[->,font=\scriptsize](m-1-2) edge node[auto] {$ g^{(\theta)} $} (m-2-2);
\end{tikzpicture}
\end{center}

Diagram lemma \ref{DLFM} implies that $\theta\ee\subseteq\ke(\F(\pi_{\theta\vv})\circ g)$ must hold.
\end{proof}
\end{prop}

\begin{hsatz}[First Diagram Lemma For $\F$-Graphs]\label{FirstDiagLemmaFGraph} Let $\G\he$, $\G\hz$ as well as $\G\hd$ be $\F$-graphs and $\phi:\G\he\twoheadrightarrow\G\hz$, $\psi:\G\he\rightarrow\G\hd$ homomorphisms, where $\phi$ is surjective. A unique  homomorphism $\gamma$, from $\G\hz$ to $\G\hd$ with $\gamma\circ\phi=\psi$, exists if and only if  $\ke\phi\subseteq\ke\psi$ holds.
\begin{proof} First, we prove the existence of $\gamma$.
\begin{align*}
\ke\phi\subseteq\ke\psi  \Longleftrightarrow & \ke\phi\vv\subseteq\ke\psi\vv\;\text{and}\;\ke\phi\ee\subseteq\ke\psi\ee\\
  \Longleftrightarrow &\exists\;\gamma\ee:E\hz\rightarrow E\hd,\; \gamma\vv:V\hz\rightarrow V\hd\;(\text{see lemma}\;\ref{DLFM})\\
  & \;\text{with}\;\gamma\ee\circ\phi\ee=\psi\ee\;\text{and}\;\gamma\vv\circ\phi\vv=\psi\vv\\
  \Longleftrightarrow & \;\exists\;\gamma:\G\hz\rightarrow\G\hd\;\text{, such that}\;\gamma\circ\phi=\psi\;\text{holds}
 \end{align*}
Hence, the map $\gamma$ exists and we have to show that it is a homomorphism.
\begin{figure}[H]
\begin{minipage}[h]{5cm}
We check:
\begin{align*}
\underline{g\hd\circ\gamma\ee}\circ\phi\ee&= g\hd\circ\psi\ee\\
                                 &=\F(\psi\vv)\circ g\he\\
                                 &=\F(\gamma\vv\circ\phi\vv)\circ g\he\\
                                 &=\F(\gamma\vv)\circ\F(\phi\vv)\circ g\he\\
                                 &=\underline{\F(\gamma\vv)\circ g\hz}\circ\phi\ee.\\
\end{align*}
\end{minipage}
\hspace{1cm}
\begin{minipage}[h]{5cm}

\begin{figure}[H]
\begin{center}
\begin{tikzpicture}[description/.style={fill=white,inner sep=2pt},>=stealth']
\matrix (m) [matrix of math nodes, row sep=1.35em,
column sep=0.6em, text height=1.5ex, text depth=0.25ex]
{ \F V\he &  & & &\F V\hz  \\
   & E\he & &E\hz & \\ 
   &  &E\hd & & \\
   &  & & & \\
   &  & \F V\hd& & \\};
\path[->>,font=\scriptsize](m-2-2) edge node[auto] {$ \phi\ee $} (m-2-4);
\path[dashed,->,font=\scriptsize](m-2-4) edge node[auto] {$ \gamma\ee $} (m-3-3);
\path[->,font=\scriptsize](m-2-2) edge node[auto,swap] {$ \psi\ee $} (m-3-3);
\path[->>,font=\scriptsize](m-1-1) edge [bend left=15] node[auto] {$ \F(\phi\vv)$} (m-1-5);
\path[->,font=\scriptsize](m-1-1) edge [bend right=15] node[auto,swap] {$ \F(\psi\vv)$} (m-5-3);
\path[dashed,->,font=\scriptsize](m-1-5) edge [bend left=15] node[auto] [bend right=60] {$ \F(\gamma\vv)$} (m-5-3);
\path[->,font=\scriptsize](m-2-2) edge node[auto,swap] {$ g\he$} (m-1-1);
\path[->,font=\scriptsize](m-2-4) edge node[auto] {$ g\hz$} (m-1-5);
\path[->,font=\scriptsize](m-3-3) edge node[description] {$ g\hd$} (m-5-3);
\end{tikzpicture}
\end{center}
\end{figure}

\end{minipage}
\end{figure}
Additionally, we know that $\phi\ee$ is surjective and thus right cancellable. 
\end{proof}
\end{hsatz}

\begin{folg}
As in lemma \ref{FirstDiagLemmaFGraph}, let $\G\he$, $\G\hz$ and $\G\hd$, together with a surjective homomorphism $\phi:\G\he\twoheadrightarrow\G\hz$, be given. If no homomorphism from $\G\hz\rightarrow\G\hd$ exists, then, for all homomorphisms $\psi:\G\he\rightarrow\G\hd$, it holds that $\ke\phi\nsubseteq\ke\psi$.
\end{folg}

\begin{hsatz}[Second Diagram Lemma For $\F$-Graphs]\label{SecondDiagLemmaFGraph} Let $\G\he$, $\G\hz$ as well as $\G\hd$ be $\F$-graphs and $\phi:\G\hz	\rightarrowtail\G\he$, $\psi:\G\hd\rightarrow\G\he$ homomorphisms, such that $\phi$ is injective. A unique  homomorphism $\gamma$, from $\G\hd$ to $\G\hz$ with $\phi\circ\gamma=\psi$, exists if and only if $\psi[\G\hd]\subseteq\phi[\G\hz]$, \textit{i.e.}, $\psi\ee{[}E\he{]}\subseteq\phi\ee{[}E\hz{]}$ and $\psi\vv{[}V\he{]}\subseteq\phi\vv{[}V\hz{]}$ holds.
\begin{proof} The proof is analogue to lemma \ref{FirstDiagLemmaFGraph} and uses lemma \ref{SecDiagLemmaSet}.
\end{proof}
\end{hsatz}

\begin{folg}
As in lemma \ref{SecondDiagLemmaFGraph}, let $\G\he$, $\G\hz$ and $\G\hd$, together with an injective homomorphism $\phi:\G\hz	\rightarrowtail\G\he$, be given. If no homomorphism from $\G\hd\rightarrow\G\hz$ exists, then, for all homomorphisms $\psi:\G\hd\rightarrow\G\he$, it holds that $\psi[\G\hd]\nsubseteq\phi[\G\hz]$.
\end{folg}

\section{Limits And Colimits In The Category Of \texorpdfstring{$\F$}{F}-Graphs}\label{ColimLimSection}
In this section, we will show that $\FGR$ is complete and cocomplete, by characterizing limits and colimits.
As colimits in $\FGR$ arise in a very natural way, we will treat them first. 
\subsection{Colimits}
In $\FGR$, each colimit is formed as a pair of the respective colimits in $\Set$.
\begin{satz}[Colimits]\label{ColimitTheorem}
The forgetful functor $\Uu:\FGR\rightarrow \Set\times\Set$ creates colimits.
\begin{proof}
Let $D:\textbf{I}\rightarrow\FGR$ be a diagram and $D\hi=(D\ee\hi,D\vv\hi,g\hi)$. Furthermore, $\Uu D=(D\ee,D\vv)$ has a colimit $(c\hi,C)_{i\in\Ob{\textbf{I}}}$ with $C=(C\ee,C\vv)$ and $c\hi=(c\hi\ee,c\hi\vv)$.
It is easy to verify that $(\F(c\hi\vv)\circ g\hi,\F C\vv)$ is a natural sink. Because of the universal property of $(c\hi\ee,C\ee)$, a map $g{_{_C}}$ with $\F(c\hi\vv)\circ g\hi=g{_{_C}} \circ c\hi\ee$ for all $i\in\Ob(\textbf{I})$ exists. We have to prove that $(C\ee,C\vv,g{_{_C}})$ is a colimit in $\FGR$. For that, let $(\tilde{c}\hi,\tilde{C})$ be a natural sink in $\FGR$. Consequently, we have a natural sink $(\tilde{c}\ee\hi,\tilde{C}\ee)$ in $\Set$ and a morphism $\delta\ee$ with $\delta\ee\circ c\hi\ee=\tilde{c}\hi\ee$, for all $i\in\Ob(\textbf{I})$. Analogously, we get a morphism $\delta\vv$ with $\delta\vv\circ c\hi\vv=\tilde{c}\hi\vv$. 

\begin{figure}[H]
\begin{minipage}[h]{5cm}
We calculate:
\begin{align*}
\underline{g_{\tilde{C}}\circ \delta\ee}\circ c\ee\hi&=g_{\tilde{C}}\circ\tilde{c}\hi\ee	\\
                                     &=\F(\tilde{c}\hi\vv)\circ g\hi\\
                                     &=\F(\delta\vv)\circ \F(c\vv\hi)\circ g\hi\\
                                     &=\underline{\F(\delta\vv)\circ g_{C}}\circ c\ee\hi.
\end{align*}
\end{minipage}
\hspace{1cm}
\begin{minipage}[h]{5cm}
\begin{center}
\begin{tikzpicture}[description/.style={fill=white,inner sep=2pt},>=stealth']
\matrix (m) [matrix of math nodes, row sep=3em,
column sep=2.5em, text height=1.5ex, text depth=0.25ex]
{ \tilde{C}\ee & D\hi\ee & C\ee\\
 \F \tilde{C}\vv & \F D\hi\vv& \F C\vv \\ };
\path[->,font=\scriptsize](m-1-2) edge node[auto] {$ \tilde{c}\ee\hi$} (m-1-1);
\path[->,font=\scriptsize](m-2-2) edge node[auto,swap] {$\F(\tilde{c}\vv\hi)$} (m-2-1);
\path[->,font=\scriptsize](m-1-3) edge[bend right=17] node[swap,auto] {$ \delta\ee  $} (m-1-1);
\path[->,font=\scriptsize](m-2-3) edge[bend left=17] node[auto] {$\F(\delta\vv)  $} (m-2-1);
\path[->,font=\scriptsize](m-1-2) edge node[auto,swap] {$ c\ee\hi $} (m-1-3);
\path[->,font=\scriptsize](m-2-2) edge node[auto] {$ \F(c\vv\hi) $} (m-2-3);
\path[->,font=\scriptsize,dashed](m-1-3) edge node[auto] {$ g{_{_C}} $} (m-2-3);
\path[->,font=\scriptsize](m-1-2) edge node[auto] {$g\hi$} (m-2-2);
\path[shorten > =2pt,->,font=\scriptsize](m-1-1) edge node[swap,auto] {$g_{\tilde{C}}$} (m-2-1);
\end{tikzpicture}
\end{center}
\end{minipage}
\end{figure}

Because $(c\ee\hi,C\ee)$ is an extremal epi-sink (see \cite{JoyCat}), we see that $(\delta\ee,\delta\vv)$ is an $\F$-graph homomorphism.
\end{proof} 
\end{satz}

Now, we have a blueprint of how to construct colimits in $\FGR$.
\begin{itemize}
	\item The coproduct $\G\he\amalg\G\hz$ has as its edge set the disjoint union $E\he\amalg E\hz$ and the vertex set is $V\he\amalg V\hz$. The structure map is defined component wise, \textit{i.e.}, $g(1,e)=g\he(e)$ and $g(2,e)=g\hz(e)$. Injections $e\he$ and $e\hz$ are given through the pairs $(e\he\ee,e\he\vv)$ and $(e\hz\ee,e\hz\vv)$.
	\item The coequalizer for two parallel homomorphisms $\phi,\psi:\G\he\rightarrow\G\hz$ is the $\F$-graph $\G/\Theta=(E/\Theta\ee,V/\Theta\vv,g^{(\Theta)})$ together with $\pi_{\Theta}=(\pi_{\Theta\ee},\pi_{\Theta\vv}):\G\hz\twoheadrightarrow\G/\Theta$, where $(\pi_{\Theta\ee},E/\Theta\ee)$ and $(\pi_{\Theta\vv},V/\Theta\vv)$ are the respective coequalizers in $\Set$ (see \cite{GummRand} for the construction of coeualizers in $\Set$). Because $\pi_{\Theta}$ is surjective, we have $g^{(\Theta)}([e]):=(\F(\pi_{\Theta\vv})\circ g\hz)(e)$.
	\item The pushout $(\phi\hP,\psi\hP,\G\hP)$ for $\phi:\G\rightarrow\G\he$ and $\psi:\G\rightarrow\G\hz$ can be formed in two steps. First, we construct the coproduct $(e\he,e\hz,\G\he\amalg\G\hz)$ and second the coequalizer of $e\he\circ\phi$ and $e\hz\circ\psi$. We get $\G\hP=((E\he\amalg E\hz)/\Theta\ee,(V\he\amalg V\hz)/\Theta\vv,g\hP)$ with $g\hP([i,e]):=(\F(\pi_{\Theta\vv})\circ g\hi)(e)$ for $i=1,2$ and $\phi\hP=\pi_{\Theta}\circ e\he$, as well as $\psi\hP=\pi_{\Theta}\circ e\hz$.
\end{itemize}

\begin{numbem}\label{gluing}
Each of the above sketched constructions can be extended to an arbitrary family of $\F$-graphs $(\G\hi)\ui$. Furthermore, in the sequel we will write $\Sigma\ui\G\hi$ for the coproduct of $(\G\hi)\ui$ and will refer to it as the sum of $(\G\hi)\ui$. Also note that the pushout enables us to define an amalgamation of graphs along a common subgraph $\G\sub$. 
\end{numbem}
\noindent
Using the sum of $\F$-graphs, we proof the following theorem, which is analogously to \cite[Theorem 3.3.2]{GU}.
\begin{satz}\label{UnionOfSubgraphs}
Let $(\G\hi\sub)_{i\in I}$ be a family of subgraphs of $\G$. Their union $\bigcup_{i\in I}\G\sub\hi$ is a subgraph of $\G$. 
\begin{proof}
For every $i\in I$, the inclusion maps $\inkl{\G^{\scriptscriptstyle{(i)}}\sub}{\G}:\G^{\scriptscriptstyle{(i)}}\sub\hookrightarrow\G$ are homomorphisms. Hence, for each $i\in I$, there is exactly one homomorphism $\sigma:\Sigma\ui\G\hi\sub\rightarrow\G$ with $\sigma\circ e^{\scriptscriptstyle{(i)}}=\inkl{\G^{\scriptscriptstyle{(i)}}\sub}{\G}$. Thus, we get $\sigma\ee\circ e^{\scriptscriptstyle{(i)}}\ee=\inkl{E^{\scriptscriptstyle{(i)}}\sub}{E}$ and $\sigma\vv\circ e\hi\vv=\inkl{V^{\scriptscriptstyle{(i)}}\sub}{V}$. The image from $\Sigma\ui\G\hi$ under $\sigma$ 

\[\sigma\ee {[}\Sigma E\hi{]}=\{\sigma\ee(i,e)\mid i\in I, e\in E\hi\}=\bigcup_{i\in I}(\sigma\ee\circ e\ee\hi)[E\hi]=\bigcup_{i\in I} E\hi,\]
\[\sigma\vv {[}\Sigma V\hi{]}=\{\sigma\vv(i,v)\mid i\in I, v\in V\hi\}=\bigcup_{i\in I}(\sigma\vv\circ e\vv\hi)[V\hi]=\bigcup_{i\in I} V\hi,\]
is a subgraph of $\G$ (see theorem \ref{hombildundteilgraph}).\qedhere 
\end{proof}
\end{satz}
\begin{defi}\label{cogensubgrdefi}
Let $\G=(E,V,g)$ be an $\F$-Graph. Because of theorem \ref{UnionOfSubgraphs}, for all subsets $E\sub\subseteq E$ and $V\sub\subseteq V$, there exists a largest subgraph whose edge set is contained in $E\sub$ and whose vertex set, in $V\sub$.  We refer to this subgraph as the  subgraph \emph{cogenerated} by $(E\sub,V\sub)$ and write $[E\sub,V\sub]\ug$. In case we have $P=(E\sub,V\sub)$, then $\hat{P}:=[E\sub,V\sub]\ug$.
\end{defi}

\subsection{Limits}
We continue with limits. Compared to colimits the situation is different, but if $\F$ preserves limits, then the limit in $\FGR$ is formed as a pair of the respective limits in $\Set$.
\begin{satz}[Limits]\label{limitpresevation}
If $\F$ preserves limits, then the forgetful functor $\Uu:\FGR\rightarrow \Set\times\Set$ creates limits.
\begin{proof}
Let $D:\textbf{I}\rightarrow\FGR$ be a diagram and $D\hi=(D\ee\hi,D\vv\hi,g\hi)$. Furthermore, $\Uu D=(D\ee,D\vv)$ has a limit $(L,l\hi)_{i\in\Ob{\textbf{I}}}$ with $L=(L\ee,L\vv)$ and $l\hi=(l\hi\ee,l\hi\vv)$.

\begin{center}
\begin{tikzpicture}[description/.style={fill=white,inner sep=2pt},>=stealth',scale=0.9]
\matrix (m) [matrix of math nodes, row sep=2.8em,
column sep=2.5em, text height=1.5ex, text depth=0.25ex]
{ L\ee & D\hi\ee \\
 \F L\vv &\F D\hi\vv \\ };
\path[->,font=\scriptsize](m-1-1) edge node[auto] {$ l\ee\hi$} (m-1-2);
\path[->,font=\scriptsize](m-2-1) edge node[auto,swap] {$\F(l\vv\hi)$} (m-2-2);
\path[dashed,->,font=\scriptsize](m-1-1) edge node[auto,swap] {$ g_{_L} $} (m-2-1);
\path[->,font=\scriptsize](m-1-2) edge node[auto] {$ g\hi $} (m-2-2);
\end{tikzpicture}
\end{center}

Because $\F$ preserves limits, $(\F L\vv,\F(l\vv\hi))$ is an $\F D\vv\hi$-limit. It is easy to check that $(L\ee,g\hi\circ l\hi\ee)$ is a natural source. Therefore, the universal property of $(\F L\vv,\F(l\vv\hi))$ gives rise to $g_{_L}\in\Mor(\Set)$, with $\F(l\vv\hi)\circ g_{_L}=g\hi\circ l\hi\ee$ for all $i\in\Ob(\textbf{I})$. The second part of the proof is analogously to the proof of \ref{ColimitTheorem}.
\end{proof}
\end{satz}

The preservation of arbitrary limits is a strong assumption. In contrast to universal coalgebras, it turns out that products in $\FGR$ exist without any condition on $\F$. %

\begin{theorem}\label{products} 
Let $(\G\hi=(E\hi,V\hi,g\hi))\ui$ be a family of $\F$-graphs. Their product is given as $\prod\G\hi=(E_{prod},\prod V\hi,g_{prod})$, where $E_{prod}$ and $g_{prod}:=\pb(\alpha)$ are defined through the following pullback square (in the sequel denoted as D1). The projection homomorphisms are $\pi\hi:=(\pi\hi\ee\circ\pb(\beta),\pi\hi\vv):\prod\G\hi\rightarrow\G\hi$.
 \begin{center}
\begin{tikzpicture}[state/.style={rectangle,rounded corners,draw=white,minimum height=2em,
           inner sep=2pt,text centered},->,>=stealth',descr/.style={fill=white,inner sep=2.5pt},scale=0.9]  
\node[state] (a){$E_{prod}$};
\node[state, right of=a,node distance=3cm,anchor=center](b){};
\node[state, above of=b,node distance=1.3cm,anchor=center](c){$\prod E\hi$};
\node[state, below of=b,node distance=1.3cm,anchor=center](d){$\F(\prod V\hi)$};
\node[state, right of=b,node distance=3cm,anchor=center](e){$\prod\F V\hi$};
\node[state, left of=a,node distance=3cm,anchor=center](f){$\tilde{E}_{prod}$};
\path[->,font=\scriptsize](a) edge node[auto] { $\pb(\beta)$ } (c);
\path[->,font=\scriptsize](a) edge node[auto,swap] { $\pb(\alpha)$ } (d);
\path[->,font=\scriptsize](c) edge node[auto] { $\alpha$ } (e);
\path[->,font=\scriptsize](d) edge node[auto,swap] { $\beta$ } (e);
\path[->,font=\scriptsize](f) edge[bend left=15] node[auto] { $\gamma$ } (c);
\path[->,font=\scriptsize](f) edge[bend right=15] node[auto,swap] { $\F(\delta\vv)\circ\tilde{g}_{prod}$ } (d);
\path[dashed,->,font=\scriptsize](f) edge node[descr] { $\exists!\,\delta\ee$ } (a);
\end{tikzpicture}
\end{center}

The maps $\alpha:=\left\langle g\hi\circ\pi\ee\hi\right\rangle\ui$ and $\beta:=\left\langle \F(\pi\vv\hi)\right\rangle\ui$ arise as mediating morphisms in the diagram below (in the sequel denoted as D2).

\begin{center}
\begin{tikzpicture}[state/.style={rectangle,rounded corners,draw=white,minimum height=2em,
           inner sep=2pt,text centered},->,>=stealth',descr/.style={fill=white,inner sep=2.5pt},transform shape,scale=0.9]  
\node[state] (a){$\F V\hi$};
\node[state, below of=a,node distance=3.2cm,anchor=center](b){$\F V\hj$};
\node[state, below of=a,node distance=1.6cm,anchor=center](c){};
\node[state, left of=c,node distance=2cm,anchor=center](d){$\prod\F V\hi$};
\node[state, left of=d,node distance=4cm,anchor=center](e){$\prod E\hi$};
\node[state, left of=d,node distance=2cm,anchor=center](f){};
\node[state, above of=f,node distance=1.6cm,anchor=center](g){$E\hi$};
\node[state, below of=f,node distance=1.6cm,anchor=center](h){$E\hj$};
\node[state, right of=d,node distance=6cm,anchor=center](i){$\F(\prod V\hi)$};

\path[->,font=\scriptsize](d) edge node[auto,swap] { $\pi\vv\hi$ } (a);
\path[->,font=\scriptsize](d) edge node[auto] { $\pi\vv\hj$ } (b);
\path[dashed,->,font=\scriptsize](e) edge node[descr] { $\exists!\,\alpha$ } (d);
\path[->,font=\scriptsize](e) edge node[auto] { $\pi\ee\hi$ } (g);
\path[->,font=\scriptsize](e) edge node[auto,swap] { $\pi\ee\hj$ } (h);
\path[->,font=\scriptsize](g) edge node[auto] { $g\hi$ } (a);
\path[->,font=\scriptsize](h) edge node[auto,swap] { $g\hj$ } (b);
\path[->,font=\scriptsize](i) edge node[auto,swap] { $\F(\pi\vv\hi)$ } (a);
\path[->,font=\scriptsize](i) edge node[auto] { $\F(\pi\vv\hj)$ } (b);
\path[dashed,->,font=\scriptsize](i) edge node[descr] { $\exists!\,\beta$ } (d);
\end{tikzpicture}
\end{center}

\begin{proof}
First, we notice that $E_{prod}$ together with $g_{prod}$ is unique up to isomorphism. We show that $\pi\hi$ is a homomorphism from $\prod\G\hi$ to $\G\hi$: $g\hi\circ\pi\ee\hi\circ\pb(\beta)\overset{\text{D2}}{=}\pi\vv\hi\circ\alpha\circ\pb(\beta)\overset{\text{D1}}{=}\pi\vv\hi\circ\beta\circ\pb(\alpha)\overset{\text{D2}}{=}\F(\pi\vv\hi)\circ g_{prod}.$	
\begin{center}
\begin{tikzpicture}[>=stealth',shorten >=2pt,shorten <=2pt,state/.style={rectangle,rounded corners,draw=white,minimum height=2em,inner sep=2pt,text centered},descr/.style={fill=white,inner sep=2.5pt},transform shape,scale=0.9]
\node[state] (a){$\prod E\hi$};
\node[state,below of=a,node distance=1.8cm,anchor=center](b){$\prod\F V\hi$};
\node[state,right of=a,node distance=3cm,anchor=center](c){$E\hi$};
\node[state,below of=c,node distance=1.8cm,anchor=center](d){$\F V\hi$};
\node[state,left of=a,node distance=3cm,anchor=center](e){$E_{prod}$};
\node[state,below of=e,node distance=1.8cm,anchor=center](f){$\F(\prod V\hi)$};
\node[state,right of=c,node distance=3cm,anchor=center](g){$\tilde{E}_{prod}$};
\node[state,right of=d,node distance=3cm,anchor=center](h){$\F(\tilde{V}_{prod})$};
\path[->,font=\scriptsize](a) edge node[auto,swap] {$ \pi\ee\hi $} (c);
\path[->,font=\scriptsize](b) edge node[auto] {$\pi\vv\hi$} (d);
\path[->,font=\scriptsize](a) edge node[auto] {$\alpha$} (b);
\path[->,font=\scriptsize](c) edge node[auto] {$g\hi$}(d);
\path[->,font=\scriptsize](e) edge node[auto] {$\pb(\beta)$}(a);
\path[->,font=\scriptsize](e) edge node[auto,swap] {$\pb(\alpha)$}(f);
\path[->,font=\scriptsize](f) edge node[auto] {$\beta$}(b);
\path[->,font=\scriptsize](g) edge node[auto] {$\tilde{g}_{prod}$}(h);
\path[->,font=\scriptsize](h) edge node[auto] {$\F(\tilde{\pi}\vv\hi)$}(d);
\path[->,font=\scriptsize](g) edge node[auto] {$\tilde{\pi}\ee\hi$}(c);
\path[->,font=\scriptsize](g) edge[bend right=17] node[descr] {$\gamma$}(a);
\path[shorten >=-2pt,shorten <=-2,->,font=\scriptsize](h) edge[bend left=25] node[auto] {$\F(\delta\vv)$}(f);
\path[->,font=\scriptsize](f) edge[bend right=17] node[descr] {$\F(\pi\hi\vv)$}(d);
\path[shorten >=-2pt,shorten <=-2,->,font=\scriptsize](g) edge[bend right=25] node[descr] {$\delta\ee$}(e);
\end{tikzpicture}
\end{center}
Next, we proof the universal property of $\prod\G\hi$. Hence, we assume the existence of some other product $\tilde{\G}_{prod}:=(\tilde{E}_{prod},\tilde{V}_{prod},\tilde{g}_{prod})$, together with projection homomorphisms $(\tilde{\pi}\ee\hi,\tilde{\pi}\vv\hi)$ to $\G\hi$.\\
Because of the universal property of $\prod E\hi$ and $\prod V\hi$, mediating morphisms $\gamma:\tilde{E}_{prod}\rightarrow \prod E\hi$ and $\delta\vv: \tilde{V}_{prod}\rightarrow \prod V\hi$ exist. It holds that $\pi\ee\hi\circ\gamma=\tilde{\pi}\ee\hi$ and that $\F(\pi\vv\hi)\circ\F(\delta\vv)=\F(\tilde{\pi}\vv\hi)$. 
We calculate: $\pi\vv\hi\circ\underline{\alpha\circ\gamma}\overset{\text{D2}}{=} g\hi\circ\pi\ee\hi\circ\gamma=g\hi\circ\tilde{\pi}\ee\hi=\F(\tilde{\pi}\vv\hi)\circ\tilde{g}_{prod}=\F(\pi\vv\hi\circ\delta\vv)\circ\tilde{g}_{prod}=\F(\pi\vv\hi)\circ\F(\delta\vv)\circ\tilde{g}_{prod}\overset{\text{D2}}{=}\pi\vv\hi\circ\underline{\beta\circ\F(\delta\vv)\circ\tilde{g}_{prod}}.$

Because the $(\pi\vv\hi)\ui$ are jointly mono, we know that the underlined equality holds. Consequently (see D1), there exists a unique $\delta\ee:\tilde{E}_{prod}\rightarrow E_{prod}$, such that $\pb(\alpha)\circ\delta\ee=\F(\delta\vv)\circ\tilde{g}_{prod}$. Therefore, $\delta:=(\delta\ee,\delta\vv)$ is a mediating $\F$-graph homomorphism form $\tilde{\G}_{prod}$ to $\prod\G\hi$. This $\delta$ is unique, as $\delta\ee$ and $\delta\vv$ are unique mediating morphisms themselves. 
\end{proof}
\end{theorem} 
From theroem \ref{products}, it follows that the edge set of the product is a canonical subset of $E\he\times E\hz\times\F(V\he\times V\hz)$.

\begin{numbsp}
In the product of two $\mathfrak{P}_{_{1,2}}$-graphs $\G$ and $\tilde{G}$, for each $e\in E$ with $g(e)=\{v,w\}$ and $\tilde{e}\in\tilde{E}$ with $\tilde{g}(e)=\{\tilde{v},\tilde{w}\}$, there are two edges in $\G\times\tilde{G}$, namely $(e,\tilde{e},\{(v,\tilde{v}),(w,\tilde{w})\})$ and $(e,\tilde{e},\{(v,\tilde{w}),(w,\tilde{v})\})$.
\end{numbsp}

\begin{numbem} We consider undirected loopless graphs. Let $\gamma:\G\he\times\G\hz\rightarrow K_n$ be a homomorphism to the complete graph with $n$ nodes, \textit{i.e.}, a proper coloring. If we can always find a (sub)graph $X$ and a homomorphism $\iota:X\rightarrow\G\he\times\G\hz$, such that for $i=1$ or $i=2$ the homomorphism $\pi\hi\circ\iota$ is surjective and $\ke(\pi\hi\circ\iota)\subseteq\ke(\gamma)$, then lemma \ref{FirstDiagLemmaFGraph} would imply Hedetniemi's conjecture. 
\end{numbem}

The next theorem about equalizers is due to \cite[Theorem 6.1]{GummRand}.
\begin{satz}\label{equalizerexistance} Let $(\phi\hi)\ui$ be a family of homomorphisms from $\G\he=(E\he,V\he,g\he)$ to $\G\hz=(E\hz,V\hz,g\hz)$. The equalizer in $\FGR$ is the $\F$-graph $\G\hEps:=[(E\hEps,V\hEps)]\ug$, for $E\hEps=\{e\in E\he\mid\forall\; i,j\in I:\; \phi\hi\ee(e)=\phi\hj\ee(e)\}$ and $V\hEps=\{v\in V\he\mid\forall\; i,j\in I:\; \phi\hi\vv(v)=\phi\hj\vv(v)\}$, together with the homomorphism $\phi\hEps:=\inkl{\G\hEps}{\G\he}$. Thus, the equalizer object in $\FGR$ is cogenerated by the equalizer in $\Set\times\Set$. 
\begin{proof}
Obviously, for all $i,j\in I$ we have $\phi\hi\circ\phi\hEps=\phi\hj\circ\phi\hEps$. We assume the existence of another $\F$-graph $\tilde{\G}\hEps$ together with $\tilde{\phi}\hEps:\tilde{\G}\hEps\rightarrow\G\he$, such that $\phi\circ\tilde{\phi}\hEps=\psi\circ\tilde{\phi}\hEps$. 
\begin{figure}[H]
\begin{center}
\begin{tikzpicture}[transform shape,scale=1.2,state/.style={rectangle,rounded corners,draw=white,minimum height=2em,
           inner sep=2pt,text centered},->,>=stealth',descr/.style={fill=white,inner sep=2.5pt}]
\node[state] (a){$\G\hEps$};
\node[state, right of=a,node distance=2cm,anchor=center](b){$\G\he$};
\node[state, right of=b,node distance=1.8cm,anchor=center](c){$\G\hz$};
\node[state, below of=a,node distance=1.3cm,anchor=center](d){$\tilde{\G}\hEps$};
\path[->,font=\tiny](b.15) edge node[auto] {$ \phi\hi$} (c.165);
\path[->,font=\tiny](b.345) edge node[swap,auto] {$ \phi\hj$} (c.195);
\path[right hook->,font=\scriptsize](a) edge node[auto] {$ \phi\hEps$} (b);
\path[shorten >=-2pt,shorten <=-2pt,->,font=\scriptsize,dashed](d) edge node[auto] {$\exists!\, \delta $} (a);
\path[->,font=\scriptsize](d) edge node[swap,auto] {$\tilde{\phi}\hEps$} (b);
\end{tikzpicture}
\end{center}
\end{figure}
The image $\tilde{\phi}\hEps[\tilde{\G}\hEps]$ is a subgraph of $\G\he$ and must be contained in $\phi\hEps[\G\hEps]$. Hence, the second diagram lemma for $\F$-graphs \ref{SecondDiagLemmaFGraph} yields a unique homomorphism $\delta:\tilde{\G}\hEps\rightarrow\G\hEps$, such that $\phi\hEps\circ\delta=\tilde{\phi}\hEps$.
\end{proof}
\end{satz}

\begin{satz}
The category $\FGR$ is complete and cocomplete. 
\begin{proof}
Follows from theorem \ref{ColimitTheorem}, \ref{products} and \ref{equalizerexistance}.
\end{proof}
\end{satz}

Now, we are able to define intersections of subgraphs and preimages of $\F$-graph homomorphisms, by means of pullbacks. 

\begin{defi}\label{intersection}
Let $(\G\hi\sub)\ui$ be a family of subgraphs of $\G$. Their \emph{intersection} is defined as the pullback of $(\inkl{\G\sub\hi}{\G}:\G\sub\hi\hookrightarrow\G)\ui$.\\ For any graph $\G=(E,V,g)$ and for all subsets $E\sub\subseteq E$ and $V\sub\subseteq V$, we define the $\F$-graph \emph{generated} by $E\sub$ and $V\sub$, \[\left\langle E\sub,V\sub\right\rangle:=\bigcap\ui\{\G\hi\sub=(E\hi\sub,V\hi\sub,g\hi\sub)\leq\G\mid E\sub\in E\hi\sub,V\sub\in V\hi\sub\},\]
\textit{i.e.}, the intersection of all subgraphs containing $E\sub$ and $V\sub$. Especially, we define the subgraph \emph{induced by one edge} $e$, namely $\G\ee=\bigcap\ui\{\G\sub=(E\hi\sub,V\hi\sub,g\sub)\leq\G\mid e\in E\sub\}$. If the type functor $\F$ preserves arbitrary intersections, then we can directly define an $\F$-graph structure on the intersection of the edge and vertex sets (see theorem \ref{limitpresevation}). Because every set endofunctor preserves nonempty finite intersections, this holds especially for all finite $\F$-graphs. 
\end{defi}
 
\begin{defi}\label{preimage}
For $\phi:\G\he\rightarrow\G\hz$ and $G\hz\sub\leq\G\hz$, the \emph{preimage} $\phi\hme[\G\hz\sub]$ is the pullback of $\inkl{\G\sub\hz}{\G\hz}$ along $\phi$. In general, this formal preimage does not have to be the preimage of $\phi\ee$ and $\phi\vv$. It follows from theorem \ref{limitpresevation} that this would be the case if $\F$ preserves pullbacks along injective maps, \textit{i.e.}, if $\F$ preserves preimages.
\end{defi}
\begin{numbem}
As in \cite{GummWeakLimit}, the preservation of certain types of limits by $\F$ can be generalized to weak preservation\footnote{In complete categories, this is equivalent to the preservation of weak limits (see \cite{GummWeakLimit}). Weak limit means that the mediating morphism is not unique.}. Especially, the weak preservation of pullbacks leads to additional structural results (as for example theorem \ref{PbGrpRelTh} or \ref{injectivemono}). In case of intersections and preimages, weak preservation by $\F$ is equivalent to preservation (see \cite{TobiSchDiss}).\\ For instance, the power set functor $\mathfrak{P}$, the identity functor and all polynomial functors weakly preserve pullbacks. Also weak limit preservation is stable under composition $\circ$, product $\times$ and sum $+$ of weak limit preserving functors (see \cite{GummWeakLimit}). 
\end{numbem}

Last, we use the union (theorem \ref{UnionOfSubgraphs}) and intersection (definition \ref{intersection}) of subgraphs, to define a complete lattice on the subgraphs of an $\F$-graph $\G$. 

\begin{prop}
The subgraphs of an $\F$-graph $\G$ define a complete bounded lattice $\Sub(\G)$ with bottom element $\G^{\scriptstyle{(\emptyset)}}=(\emptyset,\emptyset,\rho^{\scriptstyle{(\emptyset)}})$ and top element $\G$.
\begin{proof}
The subgraphs of $\G$ are partially ordered via "$\leq$". Let $(\G\hi\sub)\ui$ be subgraphs of $\G$. For the supremum and infimum operations, we define: $\bigvee\ui\G\sub\hi:=\bigcup\ui\G\hi\sub$ and $\bigwedge\ui\G\sub\hi:=\bigcap\ui\G\hi\sub$.
\end{proof}
\end{prop}

\section{Graph Relations}\label{GraphrelSection}
The concept of \enquote{bisimulation} is central for $\F$-coalgebras. Therefore, we will develop this notion for $\F$-graphs, following the presentation in \cite{GU}. Because the term \enquote{simulation} is not appropriate in the context of $\F$-graphs, we prefer to speak of \enquote{graph relations}. This construction provides a categorical foundation for the graph relations treated in \cite{LongDiss}.

\begin{defi}
Let $\G\he$ and $\G\hz$ be $\F$-graphs, and $R=(R\ee,R\vv)$ a pair of binary relations with $R\ee\subseteq E\he\times E\hz$ and $R\vv\subseteq V\he\times V\hz$. We call $R$ a \emph{graph relation} if an $\F$-graph structure $g\hr$ exist, such that the canonical projections $\pi\he=(\pi\he\ee,\pi\he\vv)$ and $\pi\hz=(\pi\hz\ee,\pi\hz\vv)$ are homomorphisms. The respective $\F$-graph will be denoted by $\G\hr=(E\hr,V\hr,g\hr):=(R\ee,R\vv,g\hr)$. If $\pi\he$ and $\pi\hz$ are surjective, we call $R$ a \emph{total graph relation}.
\begin{center}
\begin{tikzpicture}[description/.style={fill=white,inner sep=2pt},>=stealth']
\matrix (m) [matrix of math nodes, row sep=3em,
column sep=2.5em, text height=1.5ex, text depth=0.25ex]
{ E\he & R\ee & E\hz\\
 \F V\he & \F R\vv & \F V\hz \\ };

\path[->,font=\scriptsize](m-1-2) edge node[swap,auto] {$ \pi\he\ee$} (m-1-1);
\path[->,font=\scriptsize](m-1-2) edge node[auto] {$ \pi\hz\ee $} (m-1-3);
\path[->,font=\scriptsize](m-1-1) edge node[swap,auto] {$ g\he  $} (m-2-1);
\path[dashed,->,font=\scriptsize](m-1-2) edge node[auto] {$ g\hr  $} (m-2-2);
\path[->,font=\scriptsize](m-1-3) edge node[auto] {$ g\hz $} (m-2-3);
\path[->,font=\scriptsize](m-2-2) edge node[swap,auto] {$ \F(\pi\hz\vv) $} (m-2-3);
\path[->,font=\scriptsize](m-2-2) edge node[auto] {$ \F(\pi\he\vv)  $} (m-2-1);

\end{tikzpicture}
\end{center}
\end{defi}

A special case of a graph relation is the graph of a homomorphism $\phi:\G\he\rightarrow\G\hz$. We define:
\begin{itemize}
\item $G(\phi\ee):=\{(e,\phi\ee(e))\mid e\in E\he\}$ the graph of $\phi\ee$,
\item $G(\phi\vv):=\{(v,\phi\vv(v))\mid v\in V\he\}$ the graph of $\phi\vv$.
\end{itemize}
Then $G(\phi):=(G(\phi\ee),G(\phi\vv))$ is a graph relation. We even have:

\begin{satz}\label{GraphIstBisim}
A mapping $\phi:\G\he\rightarrow\G\hz$ between to $\F$-graphs is an $\F$-graph homomorphism if and only if its graph $G(\phi)$ is a graph relation between $\G\he$ and $\G\hz$.
\begin{proof}
The maps $\pi\he\ee:G(\phi\ee)\rightarrow E\he$ and $\pi\he\vv:G(\phi\vv)\rightarrow V\he$ are bijective with inverse $(\pi\he\ee)\hme$ and $(\pi\he\vv)\hme$. If $G(\phi)$ is a graph relation, then $\pi\he=(\pi\he\ee,\pi\he\vv)$ and $\pi\hz=(\pi\hz\ee,\pi\hz\vv)$ are homomorphisms. From theorem \ref{IsoCharacterization}, it will follow that $(\pi\he)\hme=((\pi\he\ee)\hme,(\pi\he\vv)\hme)$ is a homomorphism too. Hence, $\phi=\pi\hz\circ(\pi\he)\hme$ is also a homomorphism.\\
For the converse, let $\phi$ be a homomorphism. We define the structure map $g\hr:=\F((\pi\vv\he)\hme)\circ g\he\circ\pi\he\ee$
and show that $G(\phi)$ is a graph relation. For $\pi\he$ we have:
$\F(\pi\he\vv)\circ g\hr=\F(\pi\he\vv)\circ\F((\pi\he\vv)\hme)\circ g\he\circ\pi\he= g\he\circ\pi\he\ee.$
And for $\pi\hz$ we get: $\F(\pi\hz\vv)\circ g\hr=\F(\pi\hz\vv)\circ\F((\pi\he\vv)\hme)\circ g\he\circ\pi\he\ee=\F(\phi\vv)\circ g\he\circ\pi\he\ee=g\hz\circ\phi\ee\circ\pi\ee\he=g\hz\circ\pi\ee\hz.$\qedhere
\end{proof}
\end{satz}

\begin{numbsp}\label{graphofahom}
We see that the graph of a homomorphism carries a graph structure induced by the domain graph. We consider $\phi=(\phi\ee,\phi\vv)$ with $e_1\mapsto\tilde{e}_1,\,e_2\mapsto\tilde{e}_2,\,e_3\mapsto\tilde{e}_2,\,e_4\mapsto\tilde{e}_1$ and $v_1\mapsto\tilde{v}_1,\,v_2\mapsto\tilde{v}_2,\,v_3\mapsto\tilde{v}_3,\,v_4\mapsto\tilde{v}_2$.
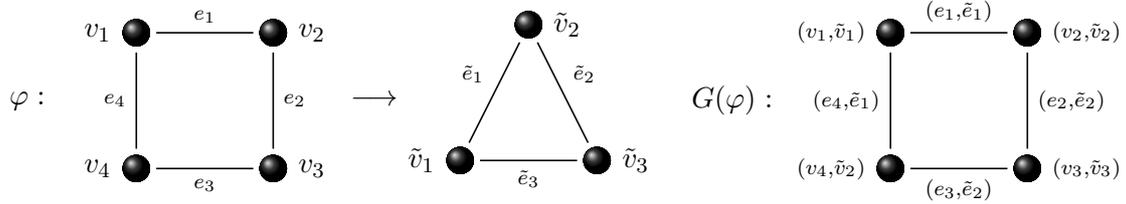
\begin{figure}[H]
\begin{minipage}[h]{0.5cm}
     $\phi:$
\end{minipage}
\hspace{0.1cm}
\begin{minipage}[h]{3.5cm}
\begin{tikzpicture}[>=stealth',shorten >=2pt , shorten <=2pt,auto,node distance=2.8cm,
                    semithick,scale=0.6,every fit/.style={ellipse,draw,inner sep=-2pt}]
   	
 	  \node[minimum size=0.3cm,circle,ball color=black,label=left:{$v_1$}](a) at (0,0){};  
 	  \node[minimum size=0.3cm,circle,ball color=black,label=right:{$v_2$}] (b) at (3,0){};
 	  \node[minimum size=0.3cm,circle,ball color=black,label=right:{$v_3$}] (c) at (3,-3){};
		\node[minimum size=0.3cm,circle,ball color=black,label=left:{$v_4$}] (d) at (0,-3){}; 
 	 
 	 \path[-,font=\scriptsize](a) edge node[auto]{$e_1$}(b);
 	 \path[-,font=\scriptsize](b) edge node[auto]{$e_2$}(c);
 	 \path[-,font=\scriptsize](c) edge node[auto]{$e_3$}(d);
 	 \path[-,font=\scriptsize](d) edge node[auto]{$e_4$}(a);
\end{tikzpicture}
\end{minipage}
\begin{minipage}[h]{0.5cm}
$\longrightarrow$
\end{minipage}
\begin{minipage}[h]{3.5cm}
\begin{tikzpicture}[>=stealth',shorten >=2pt , shorten <=2pt,auto,node distance=2.8cm,
                    semithick,scale=0.6,every fit/.style={ellipse,draw,inner sep=-2pt}]   	
 	  \node[minimum size=0.3cm,circle,ball color=black,label=left:{$\tilde{v}_1$}](a) at (0,0){};  
 	  \node[minimum size=0.3cm,circle,ball color=black,label=right:{$\tilde{v}_2$}] (b) at (1.5,3){};
 	  \node[minimum size=0.3cm,circle,ball color=black,label=right:{$\tilde{v}_3$}] (c) at (3,0){};
		 	 
 	 \path[-,font=\scriptsize](a) edge node[auto]{$\tilde{e}_1$}(b);
 	 \path[-,font=\scriptsize](b) edge node[auto]{$\tilde{e}_2$}(c);
 	 \path[-,font=\scriptsize](c) edge node[auto]{$\tilde{e}_3$}(a);
 	 \end{tikzpicture}
\end{minipage}
\hspace{0.1cm}
\begin{minipage}[h]{0.5cm}
     $G(\phi):$
 \end{minipage}
\hspace{0.5cm}
\begin{minipage}[h]{4cm}
\begin{tikzpicture}[>=stealth',shorten >=2pt , shorten <=2pt,auto,node distance=2.8cm,
                    semithick,scale=0.6,every fit/.style={ellipse,draw,inner sep=-2pt}]
   	
 	  \node[minimum size=0.3cm,circle,ball color=black,label=left:{\scriptsize ${(}v_1{,}\tilde{v}_1{)}$}](a) at (0,0){};  
 	  \node[minimum size=0.3cm,circle,ball color=black,label=right:{\scriptsize ${(}v_2{,}\tilde{v}_2{)}$}] (b) at (3,0){};
 	  \node[minimum size=0.3cm,circle,ball color=black,label=right:{\scriptsize ${(}v_3{,}\tilde{v}_3{)}$}] (c) at (3,-3){};
		\node[minimum size=0.3cm,circle,ball color=black,label=left:{\scriptsize ${(}v_4{,}\tilde{v}_2{)}$}] (d) at (0,-3){};  	 
 	 \path[-,font=\scriptsize](a) edge node[auto]{${(}e_1{,}\tilde{e}_1{)}$}(b);
 	 \path[-,font=\scriptsize](b) edge node[auto]{${(}e_2{,}\tilde{e}_2{)}$}(c);
 	 \path[-,font=\scriptsize](c) edge node[auto]{${(}e_3{,}\tilde{e}_2{)}$}(d);
 	 \path[-,font=\scriptsize](d) edge node[auto]{${(}e_4{,}\tilde{e}_1{)}$}(a);
\end{tikzpicture}
\end{minipage}
\end{figure}
\end{numbsp}

\begin{prop}
For any $\F$-graph $\G=(E,V,g)$, the diagonal $\triangle_{\scriptscriptstyle{\G}}:=(\triangle_{\scriptscriptstyle{E}},\triangle_{\scriptscriptstyle{V}})$ is a graph relation.
\begin{proof}
We consider $id:\G\rightarrow\G$. Then $G(id)=\triangle_{\scriptscriptstyle{\G}}$ is the graph of $id$ and hence a graph relation.
\end{proof}
\end{prop}

\begin{satz}\label{CharakterisierungBisimulation} Let $\G\he$ and $\G\hz$ be $\F$-graphs. For any $\F$-graph $\G$, let $\phi\he:\G\rightarrow\G\he$ and $\phi\hz:\G\rightarrow\G\hz$ be $\F$-graph homomorphisms. Then
\[(\phi\he\ee,\phi\hz\ee)[E]:=\{(\phi\he\ee(e),\phi\hz\ee(e))\mid e\in E\}\;\text{and}\]
\[(\phi\he\vv,\phi\hz\vv)[V]:=\{(\phi\he\vv(v),\phi\hz\vv(v))\mid v\in V\}\]
define a graph relation between $\G\he$ and $\G\hz$. Indeed, every graph relation arises as the image of two $\F$-graph homomorphisms.
\begin{proof}
The first claim follows from the universal property of graph products (see theorem \ref{products}).\\
For a given graph relation $R=(R\ee,R\vv)$ between $\G\he$ and $\G\hz$, there exists a graph structure $\G\hr=(E\hr,V\hr,g\hr)$, such that $\pi\he$ and $\pi\hz$ are homomorphisms. Obviously, $(\pi\he\ee,\pi\hz\ee)[E\hr]=E\hr\subseteq E\he\times E\hz$ and $(\pi\he\vv,\pi\hz\vv)[V\hr]=V\hr\subseteq V\he\times V\hz$.
\end{proof}
\end{satz}
\noindent
Depending on the type functor $\F$, pullbacks can be a source of graph relations. 
\begin{satz}\label{PbGrpRelTh} Let $\phi:\G\he\rightarrow\G\hd$ and $\psi:\G\hz\rightarrow\G\hd$ be homomorphisms. If $\F$ weakly preserves pullbacks, then the pullback of $\phi\ee,\psi\ee$ and $\phi\vv,\psi\vv$ in $\Set$ defines a graph relation between $\G\he$ and $\G\hz$. 
\begin{proof}
The proof is similar to the one of theorem \ref{limitpresevation}.
\end{proof} 
\end{satz}

\begin{folg}\label{kernelstructure} If $F$ weakly preserves kernels, then the kernel of a homomorphism $\phi:\G\he\rightarrow\G\hz$ carries a graphs structure and defines a graph relation $K(\phi)$. Furthermore, $\G\he$ is a retract of $K(\phi)$.
\end{folg}

\begin{numbsp} We consider $\phi$ from example \ref{graphofahom}.
\begin{figure}[H]
\begin{flushleft}
\begin{minipage}[h]{0.8cm}
     $\ke\phi:$ 	   
\end{minipage}
\hspace{0.3cm}
\begin{minipage}[h]{4cm}
\begin{center}

\begin{tabular}{|c|c|c|c|c|}
\hline
          &    $e_1$     &   $e_2$      & $e_3$ & $e_4$   \\ 
\hline
$e_1$     &$\times$        &             &          &  $\times$\\
\hline
$e_2$     &              &    $\times$ & $\times$ &  \\
\hline
$e_3$     &              &   $\times$  & $\times$  &  \\
\hline 
$e_4$     &$\times$    &             &          &$\times$\\
\hline                
\end{tabular}
\mdw{space}
\begin{tabular}{|c|c|c|c|c|}
\hline
          &    $v_1$     &   $v_2$      & $v_3$ & $v_4$   \\ 
\hline
$v_1$     &$\times$  &             &            & \\
\hline
$v_2$     &          &    $\times$ &            &  $\times$\\
\hline
$v_3$     &          &              &  $\times$  &  \\
\hline 
$v_4$     &          &     $\times$ &            & $\times$ \\
\hline                
\end{tabular}
\end{center}
\end{minipage}
\hspace{0.5cm}
\begin{minipage}[h]{0.5cm}
     $K(\phi):$
 \end{minipage}
\hspace{0.5cm}
\begin{minipage}[h]{4cm}
\begin{center}
\begin{tikzpicture}[scale=1.4,state/.style={rectangle,rounded corners,draw=white,minimum height=2em,
           inner sep=2pt,text centered},->,>=stealth',descr/.style={fill=white,inner sep=2.5pt}]
 	  \node[minimum size=0.3cm,circle,ball color=black,label=left:{${(}v_1{,}v_1{)}$}](a) at (0,0){};  
 	  \node[minimum size=0.3cm,circle,ball color=black,label=right:{${(}v_2{,}v_2{)}$}] (b) at (3,0){};
 	  \node[minimum size=0.3cm,circle,ball color=black,label=right:{${(}v_3{,}v_3{)}$}] (c) at (3,-3){};
		\node[minimum size=0.3cm,circle,ball color=black,label=left:{${(}v_4{,}v_4{)}$}] (d) at (0,-3){}; 
		\node[minimum size=0.3cm,circle,ball color=black,label=below:{${(}v_4{,}v_2{)}$}] (e) at (1,-2){};
		\node[minimum size=0.3cm,circle,ball color=black,label=above:{${(}v_2{,}v_4{)}$}] (f) at (2,-1){};
 	 
 	 \path[-,font=\scriptsize](a) edge node[auto]{${(}e_1{,}e_1{)}$}(b);
 	 \path[-,font=\scriptsize](b) edge node[auto]{${(}e_2{,}e_2{)}$}(c);
 	 \path[-,font=\scriptsize](c) edge node[auto]{${(}e_3{,}e_3{)}$}(d);
 	 \path[-,font=\scriptsize](d) edge node[auto]{${(}e_4{,}e_4{)}$}(a);
	 \path[-,font=\scriptsize](a) edge node[descr]{${(}e_4{,}e_1{)}$}(e);
	 \path[-,font=\scriptsize](a) edge node[descr]{${(}e_1{,}e_4{)}$}(f);
	 \path[-,font=\scriptsize](c) edge node[descr]{${(}e_3{,}e_2{)}$}(e);
	 \path[-,font=\scriptsize](c) edge node[descr]{$(e_2,e_3)$}(f);
\end{tikzpicture}
\end{center}
\end{minipage}
\end{flushleft}
\end{figure}
\end{numbsp}
                          
Similar to the proof of \ref{UnionOfSubgraphs}, we can show that graph relations are closed under union.

\begin{folg}\label{GRelGen}
If $\G\he$ and $\G\hz$ are $\F$-graphs and $R\he\subseteq E\he\times E\hz$, $R\hz\subseteq V\he\times V\hz$ a pair of relations, then there exists a largest graph relation $[R\he,R\hz]\ug=(R\ee,R\vv)$ with $R\ee\subseteq R\he$ and $R\vv\subseteq R\hz$.
\end{folg}

\begin{folg}
Between two $\F$-graphs $\G\he$ and $\G\hz$, there exists always a largest graph relation.
\end{folg}

\begin{defi}
We will denote the largest graph relation between two $\F$-graphs $\G\he$ and $\G\hz$ by $\sim\,=(\sim\ee,\sim\vv)$. Two edges $e_1\in\G\he$ and $e_2\in\G\hz$ are \emph{related} if $(e_1,e_2)\in\, \sim\ee$.

\end{defi}
With theorem \ref{CharakterisierungBisimulation}, we can easily decide whether two edges are related.

\begin{satz}\label{ununterscheidbar} Let $\G\he=(E\he,V\he,g\he)$ and $\G\hz=(E\hz,V\hz,g\hz)$ be $\F$-graphs. Two edges $e_1\in E\he$ and $e_2\in E\hz$  are related if an $\F$-graph $\G$, together with homomorphisms $\phi:\G\rightarrow\G\he$ and $\psi:\G\rightarrow\G\hz$, exists, such that $\phi\ee(e)=e_1$ and $\psi\ee(e)=e_2$, for some $e\in E$.
\end{satz}

In short: Two edges are related if they are a homomorphic image of one edge $e$.

\begin{numbsp}
Edges in undirected and directed graphs are always related. The same is true for hypergraphs. In colored graphs two edges are related, if they, together with there nodes, have the same colors.
\end{numbsp}

For an example of non related edges, we introduce the functor $(-)^k:\Set\rightarrow\Set$, which maps a set $X$ to its $k$-times product $X^k$. We define $(-)^k_l:\Set\rightarrow\Set$, which maps $X$ to $X^k$, but with the restriction that every $\overline{x}\in X^k$ has at least $l$ equal components. For instance, $(-)^3_2$ is defined as: \[(X)^3_2:=\{(x_1,x_2,x_3)\in X^3\mid x_1=x_2\vee x_1=x_3 \vee x_2=x_3\}.\]

\begin{numbsp}
We define two graphs of type $(-)^3_2$, together with edges $e$ and $\tilde{e}$. It is straight forward to show that $e$ and $\tilde{e}$ are not related. 
\begin{center}
\begin{minipage}[h]{3 cm}
    \begin{tikzpicture}[>=stealth',shorten >=2pt , shorten <=2pt,auto,node distance=2.8cm, semithick,scale=0.8]
\node (a)at(0,0)[minimum size=0.3cm,circle,ball color=black,,label=above:$v_1$]{};
\node (b)at(2,0)[minimum size=0.3cm,circle,ball color=black,,label=above:$v_2$] {}; 	   		
 	  \path[->,font=\scriptsize](a) edge [] node[auto] {$$} (b);
 	  \path[-,font=\scriptsize](a) edge [loop left] node[auto] {$$} (a);	  	
 \end{tikzpicture} 
 $g(e)=(v_1,v_1,v_2)$
\end{minipage}
\hspace{2cm}
\begin{minipage}[h]{3 cm}
   \begin{tikzpicture}[>=stealth',shorten >=2pt , shorten <=2pt,auto,node distance=2.8cm,semithick,scale=0.8]   	
\node(a)at(0,0)[minimum size=0.3cm,circle,ball color=black,label=above:$\tilde{v}_1$]{};
\node(b)at(2,0)[minimum size=0.3cm,circle,ball color=black,label=above:$\tilde{v}_2$]{}; 	   		
 \path[->,font=\scriptsize](a) edge [] node[auto] {$$} (b);
 \path[-,font=\scriptsize](b) edge [loop right] node[auto] {$$} (b);	
 \end{tikzpicture} 
 $\tilde{g}(\tilde{e})=(\tilde{v}_1,\tilde{v}_2,\tilde{v}_2)$
\end{minipage}
\end{center}
\end{numbsp}

\section{Spezial Morphisms And Isomorphism Theorems}\label{IsoEpiMonoSection}

In the first part of this section, we will characterize isomorphisms, (regular) monomorphisms and (regular) epimorphisms. This will yield a categorical formulation of the homomorphism factorization theorem \ref{hombildundteilgraph}. In the second part, we will prove the three standard isomorphism theorems for $\F$-graphs. 

\subsection{Iso-, Epi- And Monomorphisms}

\begin{satz}\label{IsoCharacterization} 
A homomorphism $\phi:\G\he\rightarrow\G\hz$ is an isomorphism if and only if it is bijective.
\begin{proof}
A bijective map in $\Set$ has an inverse. Hence, $\phi\hme:=(\phi\ee\hme,\phi\vv\hme)$ exists. We show that $\phi\hme$ is a homomorphism: $g\he\circ\phi\ee\hme=\F(\phi\vv\hme)\circ\F(\phi\vv)\circ g\he\circ\phi\ee\hme=\F(\phi\vv\hme)\circ g\hz\circ\phi\ee\circ\phi\ee\hme=\F(\phi\vv\hme)\circ g\hz.$\qedhere
\end{proof}
\end{satz}

\begin{satz}\label{EpiCharacterization} A homomorphism $\phi:\G\he\rightarrow\G\hz$ is an epimorphism if and only if $\phi$ is surjective.
\begin{proof}
Let $\phi\ee$ and $\phi\vv$ be surjective. Consequently, they are right cancellable in $\Set$. Componentwise, we have
that $(\psi\he\ee\circ\phi\ee,\psi\he\vv\circ\phi\vv)=(\psi\hz\ee\circ\phi\ee,\psi\hz\vv\circ\phi\vv).$
As $\phi\ee$ and $\phi\vv$ are right cancellable, $(\psi\he\ee,\psi\he\vv)=(\psi\hz\ee,\psi\hz\vv)$ holds.
Thus, for homomorphisms $\psi\he,\psi\hz:\G\hz\rightarrow\G\hd$ with $\psi\he\circ\phi=\psi\hz\circ\phi$, it follows that $\psi\he=\psi\hz$.\\
For the converse, let $\phi$ be an epimorphism. Hence, $(\G\hz,\id_{\G\hz})$ is the pushout of $\phi$ in $\FGR$ (see \cite{JoyCat}). Because of theorem \ref{ColimitTheorem}, $(E\hz,\id_{E\hz})$ and $(V\hz,\id_{V\hz})$ are pushouts of $\phi\ee$ and $\phi\vv$ in $\Set$. Therefore, $\phi\ee$ and $\phi\vv$ are epimorphisms in $\Set$ and thus they are surjective.
\end{proof}
\end{satz}

The characterization of regular monomorphisms is analogously to \cite[Theorem 3.4.]{GummWeakLimit}. For that, we need the following lemma. 
\begin{hsatz}\label{injsurjhom} 
Let $\G\he$, $\G\hz$ and $\G\hd$ be $\F$-graphs. Let $\phi:(E\he,V\he)\rightarrow (E\hz,V\hz)$ and $\psi:(E\hz,V\hz)\rightarrow (E\hd,V\hd)$ be maps in $\Set\times\Set$, such that $\chi:=\psi\circ\phi$ is an $\F$-graph homomorphism. The following holds:
\begin{itemize}
\item[(i)] If $\phi$ is a surjective $\F$-graph homomorphism, then $\psi$ is an $\F$-graph homomorphism.
\item[(ii)] If $\psi$ is an injective $\F$-graph homomorphism, then $\phi$ is an $\F$-graph homomorphism.
\end{itemize}
\begin{proof}
We only show $(i)$, as $(ii)$ is analogously. First, we calculate: 
\[\underline{g\hd\circ\psi\ee}\circ\phi\ee=g\hd\circ\chi\ee=\F(\chi\vv)\circ g\he=\F(\psi\vv)\circ\F(\phi\vv)\circ g\he=\underline{\F(\psi\vv)\circ g\hz} \circ \phi\ee. \]
Because $\phi\ee$ is surjective, it is right cancellable. 
\end{proof}
\end{hsatz}

\begin{satz}\label{RegularmonoCharacterization} A homomorphism $\phi:\G\he\rightarrow\G\hz$ is a regular monomorphism if and only if $\phi$ is injective.
\begin{proof} Let $\phi$ be a regular monomorphism. From theorem \ref{equalizerexistance}, it follows that $\phi$ is injective. Conversely, assume that $\phi=(\phi\ee,\phi\vv)$ is injective. In $\Set$, left inverse morphisms $\phi\ee\hme$ and $\phi\vv\hme$ exist, such that $\phi\ee\hme\circ\phi\ee=\id_{E\he}$ and $\phi\vv\hme\circ\phi\vv=\id_{V\he}$. Using lemma \ref{SecDiagLemmaSet}, we can show that $\phi\ee$ is the equalizer of $\phi\ee\circ\phi\ee\hme$ and $\id_{E\hz}$, as well as that $\phi\vv$ is the equalizer of $\phi\vv\circ\phi\vv\hme$ and $\id_{V\hz}$ in $\Set$. \\ Let $(\G\hP,\pi\he,\pi\hz)$ be the pushout of $\phi$ with itself in $\FGR$. The universal property of $\G\hP$ yields a unique map $\chi:\G\hP\rightarrow\G\hz$ in $\Set\times\Set$, such that $\chi\circ\pi\he=\phi\circ\phi\hme$ and $\chi\circ\pi\hz=(\id_{E\hz},\id_{V\hz})$.

\begin{center}
\begin{tikzpicture}[scale=1.6,state/.style={rectangle,rounded corners,draw=white,minimum height=2em,
           inner sep=2pt,text centered},->,>=stealth',descr/.style={fill=white,inner sep=2.5pt}]  
\node[state] (a){$\G\he$};
\node[state, right of=a,node distance=2.3cm,anchor=center](b){};
\node[state, above of=b,node distance=1.5cm,anchor=center](c){$\G\hz$};
\node[state, below of=b,node distance=1.5cm,anchor=center](d){$\G\hz$};
\node[state, right of=b,node distance=2.4cm,anchor=center](e){$\G\hP$};
\node[state, right of=e,node distance=2.4cm,anchor=center](f){$\G\hz$};
\node[state, left of=a,node distance=2.4cm,anchor=center](g){$\tilde{\G}$};

\path[->,font=\scriptsize](a) edge node[auto] { $\phi$ } (c);
\path[->,font=\scriptsize](a) edge node[auto,swap] { $\phi$ } (d);
\path[->,font=\scriptsize](c) edge node[auto] { $\pi\he$ } (e);
\path[->,font=\scriptsize](d) edge node[auto,swap] { $\pi\hz$ } (e);
\path[dashed,->,font=\scriptsize](e) edge node[descr] { $\chi$ } (f);
\path[->,font=\scriptsize](c) edge[bend left=15] node[auto] { $(\phi\ee\circ\phi\hme\ee{,}\phi\vv\circ\phi\hme\vv)$ } (f);
\path[->,font=\scriptsize](d) edge[bend right=15] node[auto,swap] { $(\id_{E\hz}{,}\id_{V\hz})$ } (f);
\path[dashed,->,font=\scriptsize](g) edge node[descr] { $\delta$ } (a);
\path[->,font=\scriptsize](g) edge[bend left=15] node[auto] { $\gamma$ } (c);
\path[->,font=\scriptsize](g) edge[bend right=15] node[auto,swap] { $\gamma$ } (d);

\end{tikzpicture}
\end{center}

We show that $\phi$ is the equalizer of $\pi\he$ and $\pi\hz$ in $\FGR$. For a given pair $(\tilde{\G},\gamma)$ with $\pi\he\circ\gamma$ and $\pi\hz\circ\gamma$, it follows that: $\phi\circ\phi\hme\circ\gamma=\chi\circ\pi\he\circ\gamma=\chi\circ\pi\hz\circ\gamma=(\id_{E\hz},\id_{V\hz})\circ\gamma$.

The universal property of $\phi$ as an equalizer in $\Set\times\Set$ induces a unique map $\delta:\tilde{G}\rightarrow\G\he$ with $\phi\circ\delta=\gamma$. From lemma \ref{injsurjhom}, it follows that $\delta$ is an $\F$-graph homomorphism. 
\end{proof}
\end{satz}
The characterization of monomorphisms is analogously to \cite{GU}.
\begin{satz}\label{monosatz} A homomorphism $\phi:\G\he\rightarrow\G\hz$ is a monomorphism if and only if $[\ke\phi]\ug=(\triangle_{\scriptscriptstyle{E}},\triangle_{\scriptscriptstyle{V}})$, \textit{i.e.}, the largest graph relation contained in $\ke\phi$ is the diagonal (see corollary \ref{GRelGen}).
\begin{proof}
As $[\ke\phi]\ug$ is a graph relation, the projections $\pi\he,\pi\hz:[\ke\phi]\ug\rightarrow\G\he$ are homomorphisms. Additionally, $\phi\circ\pi\he=\phi\circ\pi\hz$. If $\phi$ is mono, then $\pi\he=\pi\hz$ and hence we have that $[\ke\phi]\ug=(\triangle_{\scriptscriptstyle{E}},\triangle_{\scriptscriptstyle{V}})$.\\
For the converse, let $[\ke\phi]\ug$ be trivial. We consider homomorphisms $\psi\he,\psi\hz:\G\hd\rightarrow\G\he$, such that $\phi\circ\psi\he=\phi\circ\psi\hz$. Componentwise, we have that $(\psi\he\ee,\psi\hz\ee)[E\hd]\subseteq\ke\phi\ee$ and $(\psi\he\vv,\psi\hz\vv)[V\hd]\subseteq\ke\phi\vv$. From theorem \ref{CharakterisierungBisimulation}, we know that $(\psi\he,\psi\hz)[E\hd,V\hd]\subseteq[\ke\phi]\ug$. Because of $[\ke\phi]\ug=(\triangle_{\scriptscriptstyle{E}},\triangle_{\scriptscriptstyle{V}})$, it follows that $\psi\he=\psi\hz$. 
\end{proof}
\end{satz}
Theorem \ref{monosatz} states that a homomorphism $\phi$ is mono if and only if the largest graph relation contained in $\ke \phi$ is the diagonal. If $\phi$ is injective, then $\ke\phi=(\triangle_{\scriptscriptstyle{E}},\triangle_{\scriptscriptstyle{V}})$ and hence $\phi$ is mono. On the other hand, non injective monomorphisms can exist.

\begin{numbsp}\label{noninjectivemono}
We consider directed graphs and a non injective homomorphism $\phi$. 
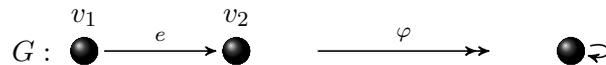
\begin{center}
\begin{tikzpicture}[>=stealth',shorten >=2pt , shorten <=2pt,auto,node distance=2.8cm, semithick,scale=0.8]
\node (a)at(0,0)[minimum size=0.3cm,circle,ball color=black,label=left:$\G:\;$,label=above:$v_1$]{};
\node (b)at(2.5,0)[minimum size=0.3cm,circle,ball color=black,label=above:$v_2$] {};
\node (c)at(8,0)[minimum size=0.3cm,circle,ball color=black] {};  	   		
 	  \path[shorten >=25pt,shorten <=25pt,->>,font=\scriptsize](b) edge node[auto] {$\phi$}(c);
 	 \path[->,font=\scriptsize](a) edge node[auto] {$e$}(b); 	 	
		\path[-,font=\scriptsize](c) edge[loop right] node[auto]{$$}(c);	  	
 \end{tikzpicture} 
\end{center}
The structure map of $\G$ is given as $g(e)=(v_1,v_2)$. If isolated nodes are forbidden, then we see that the largest graph relation contained in $\ke\phi$ is the diagonal relation. In case isolated vertices are permitted, we can include $(v_1,v_2)$ and $(v_2,v_1)$ and thus $\phi$ would be no monomorphism.
\end{numbsp}

\begin{folg}\label{injectivemono} If $\F$ weakly preserves pullbacks and isolated vertices are permitted, then every monomorphism is injective. 
\begin{proof} From corollary \ref{kernelstructure}, it follows that $[\ke\phi]\ug=\ke\phi$.
\end{proof}
\end{folg}

Analogously to \cite{TobiSchDiss}, we characterize regular epimorphisms.  

\begin{satz}
Let $\phi:\G\hz\twoheadrightarrow\G\hd$ be a surjective homomorphism. Then $\phi$ is a regular epimorphism if and only if there exists a graph relation $R$, such that $\ke\phi=\left\langle R\right\rangle$, \textit{i.e.}, $\ke\phi$ is the equivalence relation generated by $R$. 
\begin{proof}
Let $\phi$ be the coequalizer of $\psi\he,\psi\hz:\G\he\rightarrow\G\hz$ in $\FGR$. Because of theorem \ref{ColimitTheorem}, we know that $\ke\phi=\left\langle (\psi\he,\psi\hz){[}E\he,V\he{]}\right\rangle$. From theorem \ref{CharakterisierungBisimulation}, it follows that $(\psi\he,\psi\hz){[}E\he,V\he{]}$ is a graph relation.\\
Conversely, let $R$ be a graph relation on $\G\hz$ with $\ke\phi=\left\langle R\right\rangle$. For $\pi\he,\pi\hz:R\rightarrow \G\hz$, we have $\phi\circ\pi\he=\phi\circ\pi\hz$ and the universal property of $\phi$ follows from lemma \ref{FirstDiagLemmaFGraph}.
\end{proof}
\end{satz}

With the characterization of epi- and regular monomorphisms, we can rephrase the factorization theorem \ref{hombildundteilgraph} in a categorical manner. For that, we need the following lemma. 

\begin{hsatz}[Diagonal Property For $\F$-Graphs]\label{Diagonal Property} In $\FGR$, each E-M-square, where $e$ is a surjective and $m$ an injective $\F$-graph-homomorphism, has a unique diagonal $d$. 

\begin{figure}[H]
\begin{center}
\begin{tikzpicture}[description/.style={fill=white,inner sep=2pt},>=stealth']
\matrix (m) [matrix of math nodes, row sep=3em,
column sep=2.5em, text height=1.5ex, text depth=0.25ex]
{ \G\he & \G\hz \\
  \G\hd &  \G^{\scriptscriptstyle(4)} \\ };

\path[->>,font=\scriptsize](m-1-1) edge node[auto] { $e$ } (m-1-2);
\path[->,font=\scriptsize](m-1-1) edge node[auto,swap] { $h$ } (m-2-1);
\path[dashed,->,font=\scriptsize](m-1-2) edge node[description] { $d $} (m-2-1);
\path[->,font=\scriptsize](m-1-2) edge node[auto] { $f $} (m-2-2);
\path[>->,font=\scriptsize](m-2-1) edge node[auto,swap] { $m$ } (m-2-2);
\end{tikzpicture}
\end{center}
\end{figure}
\begin{proof} 
As $m$ is injective we have $\ke(e)\subseteq\ke(f\circ e)\subseteq\ke(m\circ g)\subseteq\ke(h)$. From lemma \ref{FirstDiagLemmaFGraph} the existence of a unique $\F$-graph homomorphism $d:\G\hz\rightarrow\G\hd$, with $d\circ e=h$, follows. Thus, we have $m\circ d\circ e=m\circ h=f\circ e$. Because of theorem \ref{EpiCharacterization}, $e$ is right-cancellable.
\end{proof}
\end{hsatz}

\begin{theorem}\label{FGRfactorization}
The category $\FGR$ has an epi-regular mono factorization system\footnote{Consult \cite{JoyCat} for the definition of factorization systems in categories.}.
\begin{proof}
Follows from theorem \ref{hombildundteilgraph}, \ref{EpiCharacterization} and \ref{monosatz}, and lemma \ref{Diagonal Property}. 
\end{proof}
\end{theorem}


\subsection{Isomorphism Theorems}
We will prove the three standard isomorphism theorems for $\F$-graphs. 
\begin{satz}[First Isomorphism Theorem]\label{FirstIso}
For each homomorphism $\phi:\G\he\rightarrow\G\hz$ and its kernel $\theta$, we have: $\G\he/\theta\cong\Phi[\G\he]$. If $\phi$ is surjective, $\G\he/\theta\cong\G\hz$ holds.
\begin{proof} It is easy to see that the homomorphism $\psi$, defined in the proof of theorem \ref{Faktorgraph}, is bijective. Hence, $\psi$ is an isomorphism. 
\end{proof}
\end{satz}

\begin{satz}[Second Isomorphism Theorem] Let $\G$ be an $\F$-graph and $\theta\he$, $\theta\hz$ congruence relations on $\G$ with $\theta\he\subseteq\theta\hz$. A unique homomorphism $\chi:\G/\theta\he\rightarrow\G/\theta\hz$ exist, such that $\chi\circ\pi_{\theta\he}=\pi_{\theta\hz}$. Let $\theta\hd$ be the kernel of $\chi$. It follows that $(\G/\theta\he)/\theta\hd\cong\G/\theta\hz.$
\begin{center}
\begin{tikzpicture}[shorten >=-2,shorten <=-2,scale=1.4,state/.style={rectangle,rounded corners,draw=white,minimum height=2em,inner sep=2pt,text centered},->,>=stealth',descr/.style={fill=white,inner sep=2.5pt}]  
\node[state] (a){$\G/\theta\he$};
\node[state, left of=a,node distance=3.5cm,anchor=center](b){$\G$};
\node[state, below of=a,node distance=2.3cm,anchor=center](c){$\G{/}\theta\hz$};
\node[state, right of=a,node distance=3.5cm,anchor=center](d){$(\G{/}\theta\he){/}\theta\hd$};%

\path[shorten >=2,shorten <=2,->>,font=\scriptsize](b) edge node[auto] { $\pi_{\theta\he}$ } (a);
\path[shorten <=2,->>,font=\scriptsize](b) edge node[auto,swap] {$\pi_{\theta\hz}$} (c);
\path[shorten >=1,shorten <=1,->>,font=\scriptsize](a) edge node[auto] {$\pi_{\theta\hd}$} (d);
\path[dashed,->>](a) edge node[descr] {$\chi$ } (c);
\path[>->>,font=\scriptsize](d) edge node[auto] {$$} (c);
\end{tikzpicture}
\end{center}
\begin{proof}
We define the canonical projections $\pi_{\theta\he}:\G\twoheadrightarrow\G/\theta\he$ and $\pi_{\theta\hz}:\G\twoheadrightarrow\G/\theta\hz$. Lemma \ref{FirstDiagLemmaFGraph} implies the existence of $\chi$, such that $\chi\circ\pi_{\theta\he}=\pi_{\theta\hz}$.  
Let $\pi_{\theta\hd}:\G/\theta\he\rightarrow(\G/\theta\he)/\theta\hd$ be the canonical projection. From theorem \ref{FirstIso}, we get $\G/\theta\hz=\chi[\G/\theta\he]\cong(\G/\theta\he)/\theta\hd$.
\end{proof}
\end{satz}

If $\theta\hd$ is denoted as $\theta\hz/\theta\he$, then we have $(\G/\theta\he)/(\theta\hz/\theta\he)\cong\G/\theta\hz.$

\begin{satz}[Third Isomorphism Theorem]
Let $\F$ be standard (see remark \ref{standardfunctor}), $\G=(E,V,g)$ an $\F$-graph and $\theta=(\theta\ee,\theta\vv)$ a congruence on $\G$. Furthermore, let $\G\sub=(E\sub,V\sub,g\sub)$ be a subgraph of $G$. We define $E\sub^{\theta}:=\{e\in E \mid\,\exists\, e\sub\in E\sub: (e,e\sub)\in\theta\ee \}$ and $V\sub^{\theta}:=\{v\in V \mid\,\exists\, v\sub\in V\sub: (v,v\sub)\in\theta\vv \}$. The following holds:
\begin{itemize}
\item[1.] $E\sub^{\theta}$ and $V\sub^{\theta}$ give rise to a subgraph $\G\sub^{\theta}\leq\G$,
\item[2.] $\theta\cap\G\sub\times\G\sub$ is a congruence on $\G\sub$,
\item[3.] $\G\sub/(\theta\cap\G\sub\times\G\sub)\cong\G\sub^{\theta}/\theta$.
\end{itemize}

\begin{center}
\begin{tikzpicture}[shorten >=0,shorten <=0,scale=1.4,state/.style={rectangle,rounded corners,draw=white,minimum height=2em,inner sep=2pt,text centered},->,>=stealth',descr/.style={fill=white,inner sep=2.5pt}]  

\node[state] (a){$\G\sub$};
\node[state, right of=a,node distance=3.14cm,anchor=center](b){$\G\sub^{\theta}$};
\node[state, right of=b,node distance=3.14cm,anchor=center](c){$\G$};
\node[state, right of=c,node distance=3.14cm,anchor=center](d){$\G{/}\theta$};
\path[right hook->,font=\scriptsize](a) edge node[auto] { $\inkl{\G\sub}{\G\sub^{\theta}}$ } (b);
\path[right hook->,font=\scriptsize](b) edge node[auto] {$\inkl{\G\sub^{\theta}}{\G}$ } (c);
\path[->>,font=\scriptsize](c) edge node[auto] {$\pi_{\theta}$} (d);

\end{tikzpicture}
\end{center}
\begin{proof}\mdw{Space}\\
1.: $e\in E^{\theta}\sub\Longrightarrow\exists\, e\sub\in E\sub:\,e\sub\,\theta\ee\, e \overset{\ref{ConChar}}{\Longrightarrow}\F(\pi\vv){[}g(e\sub){]}=\F(\pi\vv){[}g(e){]}\Longrightarrow g(e)\in\F V^{\theta}\sub$\\
2.: $\theta\cap\G\sub\times\G\sub=\ker(\pi_\theta\circ\inkl{\G\sub}{\G})$\\
3.: Using theorem \ref{FirstIso}, we get $\G\sub/(\theta\cap\G\sub\times\G\sub)\cong (\pi_\theta\circ\inkl{\G\sub}{\G})[\G\sub]=\pi_\theta[\G\sub^{\theta}]=\G\sub^{\theta}/\theta$.
\end{proof}
\end{satz}

\section{Cofree \texorpdfstring{$\F$}{F}-Graphs}\label{CofteeFGraphs}

Let $X=(X\ee,X\vv)$ a pair of sets. We will refer to $X\ee$ as \emph{edge colors} and to $X\vv$ as \emph{vertex colors}. For a given graph $\G$, a \emph{coloring} is a pair of maps $\gamma=(\gamma\ee,\gamma\vv):\Uu\G\rightarrow X$. Then $(\Cc(X),\varepsilon_X)$ is the \emph{cofree graph over} $X$ if for every coloring $\gamma:\Uu\G \rightarrow X$, there exists a unique $\F$-graph homomorphism $\overline{\gamma}:\G\rightarrow\Cc(X)$, such that $\varepsilon_X\circ \Uu(\bar{\gamma})=\gamma$.

\begin{figure}[H]
\begin{center}
\begin{tikzpicture}[description/.style={fill=white,inner sep=2pt},>=stealth']
\matrix (m) [matrix of math nodes, row sep=3.5em,
column sep=3em, text height=1.5ex, text depth=0.25ex] 
{\Uu\Cc(X\ee,X\vv)& & \Uu\G\\
 X=(X\ee,X\vv)& &  \\ };
\path[->,font=\scriptsize](m-1-1) edge node[auto,swap] {$ \varepsilon_X $} (m-2-1);
\path[dashed,->,font=\scriptsize](m-1-3) edge node[auto,swap] {$ \Uu(\bar{\gamma})$} (m-1-1);
\path[->,font=\scriptsize](m-1-3) edge node[auto] {$ \gamma=(\gamma\ee,\gamma\vv) $} (m-2-1);
\end{tikzpicture}
\end{center}
\end{figure}
In contrast to universal coalgebra, in $\FGR$, cofree objects exist for every $\F$.  
\begin{theorem}\label{cofreegraphcharacterization} 
The cofree graph over a color set $X=(X\ee,X\vv)$ is given through $\Cc(X)=(X\ee\times\F X\vv, X\vv, g\xx)$, with $g\xx:=\pi_{\F X\vv}:X\ee\times\F X\vv\rightarrow\F X\vv$ and $\varepsilon_X:=(\pi_{X\ee},\id_{X\vv})$. 
\begin{proof}
Let $\G=\evg$ be an $\F$-graph. Via $\gamma\ee:E\rightarrow X\ee$ and $\gamma\vv:V\rightarrow X\vv$, we define a coloring $\gamma=(\gamma\ee,\gamma\vv)$ of $\G$. There exists a homomorphism $\bar{\gamma}=(\bar{\gamma}\ee,\bar{\gamma}\vv):\G\rightarrow\Cc(X)$, where $\bar{\gamma}\ee:=\left\langle \gamma\ee, \F(\gamma\vv)\circ g\right\rangle$ is the unique mediating morphism for the product $X\ee\times\F X\vv$ and $\bar{\gamma}\vv:=\gamma\vv$. 
\begin{center}
\begin{tikzpicture}[>=stealth',shorten >=2pt,shorten <=2pt,state/.style={rectangle,rounded corners,draw=white,minimum height=2em,inner sep=2pt,text centered},descr/.style={fill=white,inner sep=2.5pt}]
\node[state] (a){$X\ee\times\F X\vv$};
\node[state,below of=a,node distance=1.7cm,anchor=center](b){$\F X\vv$};
\node[state,right of=a,node distance=2.5cm,anchor=center](c){$E$};
\node[state,below of=c,node distance=1.7cm,anchor=center](d){$\F V$};
\node[state,left of=a,node distance=2.5cm,anchor=center](e){$X\ee$};
\path[dashed,->,font=\scriptsize](c) edge node[auto] {$\bar{\gamma}\ee$} (a);
\path[->,font=\scriptsize,dashed](d) edge node[auto] {$ \F(\gamma\vv) $} (b);
\path[->,font=\scriptsize](c) edge node[auto] {$g$} (d);
\path[->,font=\scriptsize](a) edge node[auto,swap] {$\pi_{\F X\vv}$} (b);
\path[->,font=\scriptsize](a) edge node[auto] {$ \pi_{X\ee}$} (e);
\path[->,font=\scriptsize](c) edge[bend right=17] node[swap,auto] {$ \gamma\ee  $} (e);
\end{tikzpicture}
\end{center}
The homomorphism $\bar{\gamma}$ makes the diagram from the definition of cofree $\F$-graphs commute: $\varepsilon\xx\circ\Uu(\bar{\gamma})=(\pi_{X\ee},\id_{X\vv})\circ(\bar{\gamma}\ee,\bar{\gamma}\vv)=(\gamma\ee,\gamma\vv)=\gamma.$
To conclude that $\bar{\gamma}$ is unique, we assume the existence of a further $\F$-graph homomorphism $\gamma'=(\gamma\ee',\gamma\vv'):\G\rightarrow\Cc(X)$, such that $\varepsilon_X\circ\Uu(\gamma')=\gamma$ holds. Consequently, we get the equality $(\pi_{X\ee}\circ\gamma\ee',\id_V\circ\gamma\vv')=(\pi_{X\ee}\circ\bar{\gamma}\ee,\id_V\circ\bar{\gamma}\vv)$ and thus $\gamma\vv'=\bar{\gamma}\vv$. In addition to $\pi_{X\ee}\circ\gamma\ee'=\pi_{X\ee}\circ\bar{\gamma}\ee$, we compute:
$\pi_{\F X\vv}\circ\gamma\ee'=\F(\gamma'\vv)\circ g=\F(\bar{\gamma}\vv)\circ g=\pi_{\F X\vv}\circ\bar{\gamma\ee}.$
Because $\pi_{\F V}$ and $\pi_{X\ee}$ are jointly mono, it follows that $\gamma\ee'=\bar{\gamma}\ee$.
\end{proof}
\end{theorem}

\begin{numbem}\label{cofreeterminal}
From the definition of cofree graphs, it follows that the cofree graph over one vertex- and one edge color is terminal in $\FGR$. Similarly, $\Cc(X\ee,X\vv)$ is the terminal graph for the functor $\tilde{\F}(-):=X\ee\times\F((-)\times X\vv)$, \textit{i.e.}, for the category of graphs colored by $X\ee$ and $X\vv$. Hence, we can consider $\Cc(X\ee,X\vv)$ as all possible \enquote{behaviors} or \enquote{states} a graph can admit. 
\end{numbem}	

\begin{numbsp}
We consider undirected graphs of type $\mathfrak{P}_{_{1,2}}$. Let $X\ee=\{1,2,3\}$ and $X\vv=\{r,g,b\}$. The resulting cofree graph is pictured below, together with a graph coloring $\gamma$.

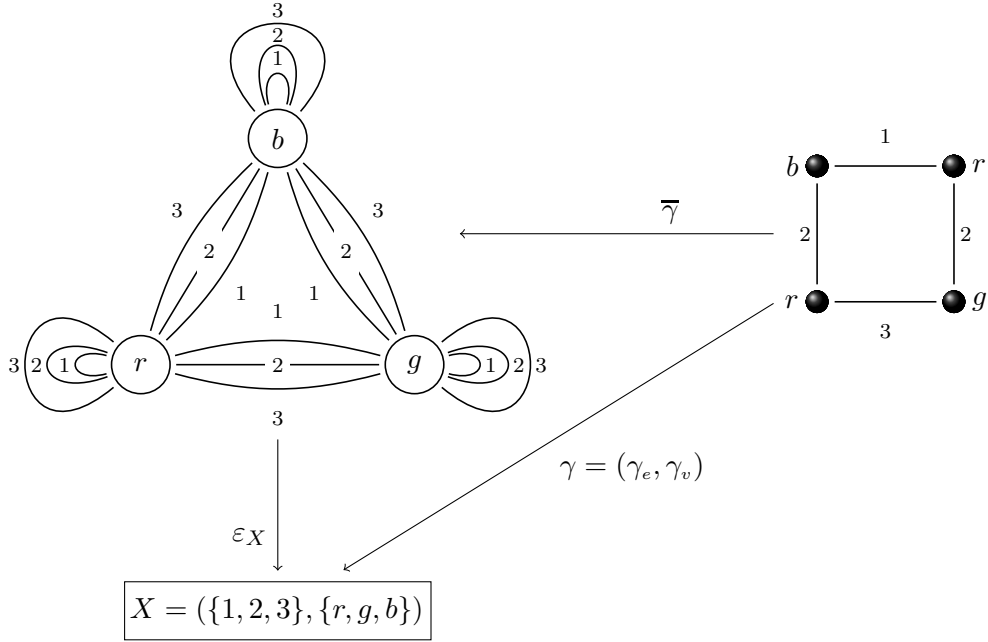
\begin{figure}[H]
\begin{center}
\begin{tikzpicture}[latticeelem/.style={rectangle,text centered,draw=white,inner sep=2pt,minimum height=2em},latticeelemblack/.style={rectangle,text centered,draw=black,inner sep=2pt,minimum height=2em},thin]

 \node[latticeelem,anchor=center](CofreeGraph)
 {
\begin{tikzpicture}[descr/.style={fill=white,inner sep=2.5pt},>=stealth',shorten >=2pt , shorten <=2pt,node distance=2.8cm,semithick,scale=0.6,every fit/.style={ellipse,draw,inner sep=-2pt},every loop/.style={}]
   	
 	  \node[draw,circle] (a) at (0,0){$r$};  
 	  \node[draw,circle] (b) at (6,0){$g$};
 	  \node[draw,circle] (c) at (3,5){$b$};	
 	
 	 \path[-,font=\scriptsize](a) edge node[descr] {$2$}(b);
 	 \path[-,font=\scriptsize](b) edge node[descr] {$2$}(c);
 	 \path[-,font=\scriptsize](c) edge node[descr] {$2$}(a);
	 \path[-,font=\scriptsize](a) edge[bend left=15] node[auto]{$1$}(b);
 	 \path[-,font=\scriptsize](b) edge[bend left=15] node[auto]{$1$}(c);
 	 \path[-,font=\scriptsize](c) edge[bend left=15] node[auto]{$1$}(a);
	 \path[-,font=\scriptsize](a) edge[bend right=15] node[auto,swap]{$3$}(b);
 	 \path[-,font=\scriptsize](b) edge[bend right=15] node[auto,swap]{$3$}(c);
 	 \path[-,font=\scriptsize](c) edge[bend right=15] node[auto,swap]{$3$}(a);
	
	 \path[-,font=\scriptsize](a) edge[loop left] node[auto]{$1$}(a);
	\path[-,font=\scriptsize](a) edge[loop left,in=200,out=160,min distance=20mm] node[auto,swap]{$2$}(a);
	\path[-,font=\scriptsize](a) edge[loop left,in=220,out=140,min distance=35mm] node[auto,swap]{$3$}(a);
	
	 \path[-,font=\scriptsize](b) edge[loop right] node[auto]{$1$} (b);
	\path[-,font=\scriptsize](b) edge[in=340,out=20,min distance=20mm] node[auto]{$2$} (b);
	\path[-,font=\scriptsize](b) edge[in=320,out=40,min distance=35mm] node[auto]{$3$} (b);
	
	 \path[-,font=\scriptsize](c) edge[loop above] node[above=-1.8mm]{$1$}(c);
	\path[-,font=\scriptsize](c) edge[in=110,out=70,min distance=20mm] node[above=-2.5mm]{$2$} (c);
	\path[-,font=\scriptsize](c) edge[in=130,out=50,min distance=35mm] node[above=-2mm]{$3$} (c);

\end{tikzpicture}
  };  
\node[latticeelem,right of=CofreeGraph,node distance=8cm,anchor=center] (TestGraph)
 {

\begin{tikzpicture}[>=stealth',shorten >=2pt , shorten <=2pt,auto,node distance=2.8cm,
                    semithick,scale=0.6,every fit/.style={ellipse,draw,inner sep=-2pt}]
   	
 	  \node[minimum size=0.3cm,circle,ball color=black,label=left:{$b$}](a) at (0,0){};  
 	  \node[minimum size=0.3cm,circle,ball color=black,label=right:{$r$}] (b) at (3,0){};
 	  \node[minimum size=0.3cm,circle,ball color=black,label=right:{$g$}] (c) at (3,-3){};
		\node[minimum size=0.3cm,circle,ball color=black,label=left:{$r$}] (d) at (0,-3){}; 
 	 
 	 \path[-,font=\scriptsize](a) edge node[auto]{$1$}(b);
 	 \path[-,font=\scriptsize](b) edge node[auto]{$2$}(c);
 	 \path[-,font=\scriptsize](c) edge node[auto]{$3$}(d);
 	 \path[-,font=\scriptsize](d) edge node[auto]{$2$}(a);
	\path[-,font=\scriptsize](a) edge[white, loop above] (a);
\end{tikzpicture}  
  };
\node[latticeelemblack,below of=CofreeGraph,node distance=5cm,anchor=center] (ColorSet)
 {
   $X=(\{1,2,3\},\{r,g,b\})$
   
  };

\path[->,shorten >=8pt](TestGraph)   edge  node[auto] {$\gamma=(\gamma\ee,\gamma\vv)$}(ColorSet);
\path[->,shorten >=4pt,shorten <=-20pt](CofreeGraph) edge  node[auto,swap] {$\varepsilon_X$} (ColorSet);
\path[->,shorten >=-40pt](TestGraph)   edge  node[auto,swap] {$\overline{\gamma}$}(CofreeGraph);
  
\end{tikzpicture}
\end{center}
\end{figure}
\end{numbsp}

The cofree property of $\Cc(X)$ implies that $\Cc:\Set\times\Set\rightarrow\FGR$ defines a right-adjoint functor $\Uu\dashv\Cc$, where $\varepsilon_X:\Uu\Cc(X)\rightarrow X$ is the adjunctions counit. We refer to $\Cc$ as the \emph{cofree graph functor}. The next lemma shows how $\Cc$ acts on morphisms. 

\begin{hsatz}\label{CofreeFunctorMorphisms} Let $\gamma=(\gamma\ee,\gamma\vv):(X\ee\he,X\vv\he)\rightarrow(X\ee\hz,X\vv\hz)$ be a morphism in $\Set\times\Set$ and $\Cc:\Set\times\Set\rightarrow\FGR$ the cofree graph functor. Then $\Cc(\gamma)$ is given as $\Cc(\gamma\ee,\gamma\vv):=(\alpha,\gamma\vv)$, where $\alpha:=\langle\gamma\ee\circ\pi_{X\ee\he},\F(\gamma\vv)\circ\pi_{\F X\he\vv}\rangle$ is the unique mediating morphism in the following diagram.
\begin{center}
\begin{tikzpicture}[scale=1.0,state/.style={rectangle,rounded corners,draw=white,minimum height=2em,inner sep=2pt,text centered},->,>=stealth',descr/.style={fill=white,inner sep=2.5pt}]  
\node[state] (a){$X\ee\hz$};
\node[state, below of=a,node distance=3cm,anchor=center](b){$\F X\vv\hz$};
\node[state, below of=a,node distance=1.5cm,anchor=center](c){};
\node[state, left of=c,node distance=2cm,anchor=center](d){$X\ee\hz\times\F X\vv\hz$};
\node[state, left of=d,node distance=4cm,anchor=center](e){$X\ee\he\times\F X\vv\he$};
\node[state, left of=d,node distance=2cm,anchor=center](f){};
\node[state, above of=f,node distance=1.5cm,anchor=center](g){$X\ee\he$};
\node[state, below of=f,node distance=1.5cm,anchor=center](h){$\F X\vv\he$};

\path[->,font=\scriptsize](d) edge node[auto,swap] { $\pi_{X\ee\hz}$ } (a);
\path[->,font=\scriptsize](d) edge node[auto] { $\pi_{\F X\vv\hz}$ } (b);
\path[dashed,->,font=\scriptsize](e) edge node[descr] { $\alpha$ } (d);
\path[->,font=\scriptsize](e) edge node[auto] { $\pi_{X\ee\he}$ } (g);
\path[->,font=\scriptsize](e) edge node[auto,swap] { $\pi_{\F X\vv\he}$ } (h);
\path[->,font=\scriptsize](g) edge node[auto] { $\gamma\ee$ } (a);
\path[->,font=\scriptsize](h) edge node[auto,swap] { $\F(\gamma\vv)$ } (b);
\end{tikzpicture}
\end{center}
\begin{proof}
By definition $\pi_{\F X\hz\vv}\circ\alpha=\F(\gamma\vv)\circ\pi_{\F X\vv\he}$ holds and thus $(\alpha,\gamma\vv):\Cc(X\he\ee,X\he\vv)\rightarrow\Cc(X\hz\ee,X\hz\vv)$ is a homomorphism. Because $\alpha$ is defined through a functorial construction, the conditions $\Cc(\id_{X\ee},\id_{X\vv})=\id_{\Cc(X\ee,X\vv)}$ and $\Cc(\gamma\hz\circ\gamma\he)=\Cc(\gamma\hz)\circ\Cc(\gamma\he)$ follow. 
\end{proof}
\end{hsatz}
It is a property of cofree graphs that every $F$-graph $\G$ is a subgraph of some cofree graph.
\begin{hsatz}\label{unitembedding}
Let $\G=\evg$ be an $\F$-graph. It is a subgraph of $\Cc\Uu\G$, \textit{i.e.}, the cofree graph over $(E,V)$.
\begin{proof}
Let $\beta:=\left\langle\id_E,\F(\id_V)\circ g\right\rangle$ be the mediating morphism in the following diagram. 
\begin{center}
\begin{tikzpicture}[scale=1.6,state/.style={rectangle,rounded corners,draw=white,minimum height=2em,
           inner sep=2pt,text centered},->,>=stealth',descr/.style={fill=white,inner sep=2.5pt}]  
\node[state] (a){$E$};
\node[state, below of=a,node distance=2cm,anchor=center](b){$\F V$};
\node[state, below of=a,node distance=1cm,anchor=center](c){};
\node[state, left of=c,node distance=2cm,anchor=center](d){$E\times\F V$};
\node[state, left of=d,node distance=2cm,anchor=center](e){$E$};

\path[->,font=\scriptsize](d) edge node[auto,swap] { $\pi_E$ } (a);
\path[->,font=\scriptsize](d) edge node[auto] { $\pi_{\F V}$ } (b);
\path[dashed,->,font=\scriptsize](e) edge node[descr] { $\beta$ } (d);
\path[->,font=\scriptsize](e) edge[bend left=15] node[auto] { $\id_E$ } (a);
\path[->,font=\scriptsize](e) edge[bend right=15] node[auto,swap] { $\F(\id_V)\circ g$ } (b);
\end{tikzpicture}
\end{center}
It holds that $\id_E=\pi_E\circ\beta$. Consequently, $\beta$ has a left inverse and is therefore injective. It follows from theorem \ref{RegularmonoCharacterization} that $\G\leq\Cc\Uu\G$, via $\eta_{\G}:=(\beta,\id_V):\G\rightarrowtail\Cc(E,V)=\Cc\Uu\G$.
\end{proof}
\end{hsatz}

\begin{numbem}
As $\Cc$ is a right-adjoint functor, it preserves products. Hence, $\Cc(X\he\ee,X\he\vv)\times\Cc(X\hz\ee,X\hz\vv)=\Cc(X\he\ee\times X\hz\ee,X\he\vv\times X\hz\vv)$. From this observation, it follows that the product (see theorem \ref{products}), of two arbitrary $\F$-graphs $\G\he=(E\he,V\he,g\he)$ and $\G\hz=(E\hz,V\hz,g\hz)$, is a subgraph of $\Cc(E\he,V\he)\times\Cc(E\hz,V\hz)$. Following the approach of \cite{ProdCoalgGummSchr}, it is possible to define graph products as subgraphs of cofree graphs with a certain maximality condition.
\end{numbem}

We saw that the functors $\Uu\dashv\Cc$ define an adjunction and that $\varepsilon$ is the respective counit. Let $X$ be a set of colors. From the duality of adjoint situations, we know that a map $\eta_\G:\G\rightarrow\Cc\Uu\G$ must exist, such that for every set of colors $X=(X\ee,X\vv)$ and every homomorphism $\bar{\gamma}:\G\rightarrow\Cc X$, a unique coloring $\gamma:\Uu\G\rightarrow X$ is induced. For this coloring $\Cc(\gamma)\circ\eta_\G=\bar{\gamma}$ holds.

\begin{satz}\label{UnitCofreeGraph} The unit of $\,\Uu\dashv\Cc$ is the homomorphism $\eta_{\G}=(\beta,\id_V):\G\rightarrowtail\Cc\Uu\G$, defined in lemma \ref{unitembedding}. 
\begin{proof}
Let $X$ be a set of colors and $\bar{\gamma}=(\bar{\gamma}\ee,\bar{\gamma}\vv):\G\rightarrow\C X$ an $\F$-graph homomorphism. Therefore, $\F(\bar{\gamma}\vv)\circ g=\pi_{\F X\vv}\circ\bar{\gamma}\ee
$ holds. We claim that the unique coloring $\gamma=(\gamma\ee,\gamma\vv)$ is given via $\gamma\ee:=\pi_{X\ee}\circ\bar{\gamma}\ee$ and $\gamma\vv=\bar{\gamma}\vv$.\\
Let $\alpha:=\langle\pi_{X\ee}\circ\bar{\gamma}\ee\circ\pi_{E},\F(\bar{\gamma}\vv)\circ\pi_{\F V}\rangle$ be defined as in lemma \ref{CofreeFunctorMorphisms}. We calculate: $\Cc(\gamma)\circ\eta_\G=(\alpha,\gamma\vv)\circ(\beta,\id_V)=(\alpha\circ\beta,\bar{\gamma}\vv).$
For $\alpha\circ\beta$, it holds that: $\pi_{X\ee}\circ\alpha\circ\beta=\pi_{X\ee}\circ\bar{\gamma}\ee\circ\pi_E\circ\beta=\pi_{X\ee}\circ\bar{\gamma}\ee.$
Furthermore, we calculate: $\pi_{\F X\vv}\circ\alpha\circ\beta=\F(\bar{\gamma}\vv)\circ\pi_{\F V}\circ\beta=\F(\bar{\gamma}\vv)\circ g=\pi_{\F X\vv}\circ\bar{\gamma}\ee.$
Because $\pi_{X\ee}$ and $\pi_{\F X\vv}$ are jointly mono, it follows that $\alpha\circ\beta=\bar{\gamma}\ee$.\\ 
To proof the uniqueness of $\gamma$, we assume that $\gamma'$, such that $\Cc(\gamma')\circ\eta_\G=\bar{\gamma}$, exists. As in lemma \ref{CofreeFunctorMorphisms}, let $\alpha':=\langle\gamma\ee'\circ\pi_E,\F(\gamma\vv')\circ\pi_{\F V}\rangle$. Because of $(\alpha',\gamma\vv')\circ(\beta,\id_V)=(\alpha,\gamma\vv)\circ(\beta,\id_V)$, we know that $\gamma'\vv=\gamma\vv$. Additionally, $\alpha'\circ\beta=\alpha\circ\beta$ holds. At last, we show the equality
$\gamma\ee'=\gamma\ee'\circ\pi_E\circ\beta=\pi_{X\ee}\circ\alpha'\circ\beta=\pi_{X\ee}\circ\alpha\circ\beta=\gamma\ee\circ\pi_E\circ\beta=\gamma\ee.$
\end{proof}
\end{satz}

\begin{numbem}
If we want to consider cofree graphs over vertex colors only, there are two ways to do this. First, we could assign one similar color to all edges and form the cofree graph over $X=(\{\ast\},X\vv)$. Another possibility would be to define the \emph{underlying vertex set functor} $\forg{-}{V}:\FGR\rightarrow\Set$, via: $\forg{\G}{V}=V\;\text{and}\;\forg{\phi}{V}=\phi\vv.$
Next, we prove the existence of a right adjoint functor $\Cc\vv$. The latter approach is described in \cite{DeCatGraph}. There it is shown that for graphs of type $\mathfrak{P}_{_{1,2}}$, the cofree graph over $X\vv$ is the complete graph with $|X\vv|$ vertices. Each vertrex represents one color of $X\vv$. 

In a similar way, we define the \emph{underlying edge set functor} $\forg{-}{E}:\FGR\rightarrow\Set$:
\[\forg{\G}{E}=E\;\text{and}\;\forg{\phi}{E}=\phi\ee.\]
There exists a right adjoint functor $\Cc\ee$. Again, an identical construction is obtained if we consider the cofree graph over $X=(X\ee,\{\ast\})$. Hence, for graphs of type $\mathfrak{P}_{_{1,2}}$, the cofree graph over $X\ee$ has one vertex and one loop for each edge color.
\end{numbem}

Obviously, there is a close relationship between $\Uu$, $\forg{-}{E}$ and $\forg{-}{V}$, namely the tupeling $\Uu\G=(\forg{\G}{E},\forg{\G}{V})$. 
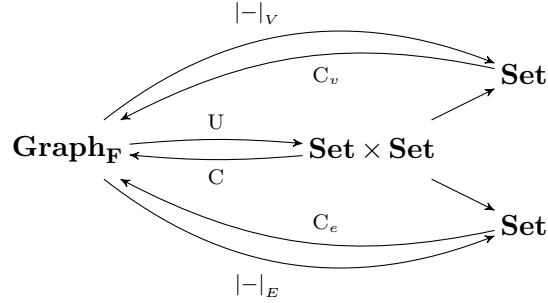
\begin{figure}[H]
\begin{center}
\begin{tikzpicture}[scale=1.6,state/.style={rectangle,rounded corners,draw=white,minimum height=2em,
           inner sep=2pt,text centered},->,>=stealth',descr/.style={fill=white,inner sep=2.5pt}]  
\node[state] (a){$\Set$};
\node[state, below of=a,node distance=2cm,anchor=center](b){$\Set$};
\node[state, below of=a,node distance=1cm,anchor=center](c){};
\node[state, left of=c,node distance=2cm,anchor=center](d){$\Set\times\Set$};
\node[state, left of=d,node distance=4cm,anchor=center](e){$\FGR$};

\path[->,font=\scriptsize](d) edge node[auto,swap] { $$ } (a);
\path[->,font=\scriptsize](d) edge node[auto] { $$ } (b);
\path[->,font=\scriptsize](e) edge[bend left=5] node[auto] { $\Uu$ } (d);
\path[->,font=\scriptsize](d) edge[bend left=5] node[auto] { $\Cc$ } (e);

\path[->,font=\scriptsize](e) edge[bend left=30] node[auto] { $\forg{-}{V}$ } (a);
\path[->,font=\scriptsize](a) edge[bend right=20] node[auto] { $\Cc\vv$ } (e);
\path[->,font=\scriptsize](e) edge[bend right=30] node[auto,swap] { $\forg{-}{E}$ } (b);
\path[->,font=\scriptsize](b) edge[bend left=20] node[auto,swap] { $\Cc\ee$ } (e);
\end{tikzpicture}
\end{center}
\end{figure}

A cofree graph can always be decomposed into a product of two graphs, induced by $X\ee$ and $X\vv$, via $\Cc\ee$ and $\Cc\vv$.
\begin{prop} Let $X\ee$ be a set of edge colors and $X\vv$ a set of vertex colors. It holds that $\Cc(X\ee,X\vv)\cong\Cc\ee(X\ee)\times\Cc\vv(X\vv)$. 
\begin{proof} $\Cc(X\ee,X\vv)\cong\Cc(X\ee\times\{\ast\},\{\ast\}\times X\vv)\cong\Cc(X\ee,\{\ast\})\times\Cc(\{\ast\},X\vv)\cong\Cc\ee(X\ee)\times\Cc\vv(X\vv).$\qedhere
\end{proof}
\end{prop}

Next, we proof an important property of cofree graphs\footnote{The proof is similar to the on in \cite{GummRand}.}.

\begin{lemma}\label{ExtensionCoFree} Let $\G\sub\leq\G$ for any $\F$-graph $\G$. Every homomorphism $\phi:\G\sub\rightarrow\Cc X $ can be extended to $\psi:\G\rightarrow\Cc X$, \textit{i.e.}, $\phi=\psi\;\circ\inkl{\G\sub}{\G}$.
\begin{proof} We define a coloring of $\G\sub$ via $\varepsilon_X\circ\Uu(\phi)$. In $\Set\times\Set$, it can be extended to a coloring $\gamma:\Uu\G\rightarrow X$. Thus, we get the homomorphic extension $\overline{\gamma}:\G\rightarrow\Cc X $. We define $\psi:=\overline{\gamma}$ and calculate: $\varepsilon_X\circ\Uu(\psi\circ\inkl{\G\sub}{\G})=\varepsilon_X\circ\Uu(\psi)\circ\Uu(\inkl{\G\sub}{\G})=\gamma\circ\Uu(\inkl{\G\sub}{\G})=\varepsilon_X\circ\Uu(\phi)$. From the uniqueness of $\overline{\gamma}$, it follows that $\psi\;\circ\inkl{\G\sub}{\G}=\phi$.
\begin{figure}[h]
\begin{center}
\begin{tikzpicture}[>=stealth',shorten >=2pt,shorten <=2pt,state/.style={rectangle,rounded corners,draw=white,minimum height=2em,inner sep=2pt,text centered},->,>=stealth',descr/.style={fill=white,inner sep=2.5pt}]
 	      	
 	  \node[state] (a){$\Uu\G$};
 	  \node[state,below of=a,node distance=1.7cm,anchor=center](b){$\Uu\G\sub$};
 	  \node[state,right of=a,node distance=2.5cm,anchor=center](c){$X$};
 	  \node[state,below of=c,node distance=1.7cm,anchor=center](d){$\Uu\Cc X$};
 	    	    
 	  \path[right hook->,shorten >=0pt,shorten <=0pt,font=\scriptsize](b) edge  node[auto] {$\Uu(\iota)$} (a); 
 	  \path[->,font=\scriptsize](a) edge  node[auto] {$\gamma$} (c);
 	  \path[->,font=\scriptsize](b) edge  node[auto,swap]  {$\Uu(\phi)$}(d); 
 	  \path[shorten >=0pt,shorten <=0pt,->,font=\scriptsize](d) edge  node[auto,swap] {$\varepsilon_X$} (c);
		\path[dashed,->,font=\scriptsize](a) edge  node[descr] {$\Uu(\psi)$} (d);
\end{tikzpicture}
\end{center}
\end{figure}
\end{proof}
\end{lemma}
Lemma \ref{ExtensionCoFree} states that cofree Graphs are regular-injective objects (see \cite{JoyCat}). 
The following theorem generalizes this fact.
\begin{theorem}\label{RegInjChar}
$\G$ is a regular-injective $\F$-graph if and only if it is a retract of some cofree graph. 
\begin{proof}
Let $\G$ be regular-injective. From lemma \ref{unitembedding}, we know that $\G$ is a subgraph of $\Cc\Uu\G$. Consequently, there exists an extension $\phi$ of $\id_{\G}:\G\rightarrow\G$, such that $\phi\circ\iota_\G=\id_\G$.\\
For the converse, let $\G$ be a retract of the cofree graph $\Cc X$. Hence, there are $\F$-graph homomorphisms $\lambda:\G\rightarrow\C X$ and $\pi:\C X\rightarrow\G$, with $\pi\circ\lambda=\id_\G$. Next, we consider $\tilde{G}\sub\leq\tilde{\G}$ together with $\phi:\tilde{\G}\sub\rightarrow\G$. We have to construct $\psi:\tilde{\G}\rightarrow\G$, such that $\psi\circ\inkl{\tilde{\G}\sub}{\tilde{\G}}=\phi$.
\begin{center}
\begin{tikzpicture}[shorten >=0,shorten <=0,scale=1.4,state/.style={rectangle,rounded corners,draw=white,minimum height=2em,inner sep=2pt,text centered},->,>=stealth',descr/.style={fill=white,inner sep=2.5pt}]  

\node[state] (a){$\tilde{\G}\sub$};
\node[state, right of=a,node distance=2.3cm,anchor=center](b){$\tilde{G}$};
\node[state, below of=b,node distance=1.7cm,anchor=center](c){$\G$};
\node[state, right of=c,node distance=2.3cm,anchor=center](d){$\Cc X$};
\path[>->,font=\scriptsize](a) edge node[auto] { $\inkl{\tilde{\G}\sub}{\tilde{\G}}$ } (b);
\path[->,font=\scriptsize](a) edge node[auto,swap] { $\phi$ } (c);
\path[dashed,->,font=\scriptsize](b) edge node[auto] {$\psi$} (c);
\path[->,font=\scriptsize](b) edge node[auto] {$\tilde{\phi}$ } (d);
\path[->,font=\scriptsize](c) edge[bend left=15] node[auto] {$\lambda$ } (d);
\path[->,font=\scriptsize](d) edge[bend left=15] node[auto] {$\pi$ } (c);
\end{tikzpicture}
\end{center}

Because $\Cc X$ is regular-injective, $\tilde{\phi}:\tilde{\G}\rightarrow\Cc X$ with $\tilde{\phi}\circ\inkl{\tilde{\G}\sub}{\tilde{\G}}=\lambda\circ\phi$ exists. Finally, we define $\psi:=\pi\circ\tilde{\phi}$ and calculate: $\psi\circ\inkl{\tilde{\G}\sub}{\tilde{\G}}=\pi\circ\tilde{\phi}\circ\inkl{\tilde{\G}\sub}{\tilde{\G}}=\pi\circ\lambda\circ\phi=\phi.$\qedhere
\end{proof}
\end{theorem}

\begin{numbem}
For the existence of free $\F$-graphs over $(X\ee,X\vv)$, a left-adjoint to the forgetful functor $\Uu:\FGR \rightarrow \Set\times\Set$ would be necessary. According to the special adjoint functor theorem, the functor $\Uu$ would have to preserve limits. Generally, this is not the case. Nevertheless, free graphs for the underlying vertex set functor $\forg{-}{V}$ can be constructed. For graphs of type $\mathfrak{P}_{_{1,2}}$, the free graph over $X\vv$ is the graph with empty edge set and $X\vv$ as its vertex set (see \cite{DeCatGraph}).
\end{numbem}

\begin{numbem}
The adjunction $\Uu\dashv\Cc$, with unit $\eta$ and counit $\varepsilon$, gives rise to a monad $\left\langle\Cc\Uu,\eta,\Cc\varepsilon_{\Uu}\right\rangle$ and a comonad $\left\langle\Uu\Cc,\varepsilon,\Uu\eta_{\Cc}\right\rangle$. For instance, the functor $\Cc\Uu:\FGR\rightarrow\FGR$ maps a graph $\G$ to the \enquote{hull} induced by $(E,V)$. 
\end{numbem}
\section{Conjunct Sums And Graph Transformations}\label{CSumsGrTrafoSection}
In the first part of this section, we will proof a dual statement to the fact that every algebra is a subdirect product of subdirectly irreducible algebras. In the second part, we will consider functors between categories of graphs.

\subsection{Conjunct Sums}
The following presentation is based on \cite[Section.7.2]{GummRand}.
\begin{defi}
Let $(\G\hi)\ui$ be a family of $\F$-graphs. A graph $\G=(E,V,g)$ is called a \emph{conjunct sum} of the $\G\hi=(E\hi,V\hi,g\hi)$ if for every $e\in E$ there is an $i\in I$ and an injective homomorphism $\phi\hi:\G\hi\rightarrowtail\G$ with $e\in\phi\ee\hi[E\hi]$. A \emph{conjunct representation of $\G$} is a family $(\phi\hk:\G\hk\rightarrowtail\G)_{k\in K\subseteq I}$, such that $\bigcup_{k\in K}\phi\hk[E\hk,V\hk]=(E,V)$.
\end{defi}

Thus, if $\G=(E,V,g)$ is a conjunct sum of the $\G\hi$, then for every $e\in E$ a subgraph $\G\sub\leq\G$  exists, which contains $e$ and is isomorphic to one of the $\G\hi$. Therefore, we can consider the $\G\hi$ as (possibly overlapping) building blocks from which $\G$ is constructed. 

Each graph has a trivial conjunct representation $\id:\G\rightarrow\G$. If there is no other representation, we call this graph conjunctly irreducible. More precisely:

\begin{defi} A graph $\G$ is called \emph{conjunctly irreducible} if in each conjunct representation $(\phi\hi:\G\hi\rightarrowtail\G)\ui$ of $\G$, one of the $\phi\hi$ is an isomorphism. 
\end{defi}

\begin{defi}
An $\F$-graph $\G=(E,V,g)$ is called \emph{one-generated} if there is an $e\in E$, such that $\G$ is the only subgraph of $\G$ containing $e$. 
\[e\in\G\sub\leq\G\Longrightarrow \G\sub=\G\]
\end{defi}

It is easy to show that the from one edge induced subgraphs (see definition \ref{intersection}) are one-generated.

\begin{prop}
A graph is conjunctly irreducible if and only if it is one-generated. 
\begin{proof}
"$\Rightarrow$": We assume that $\G=\evr$ is conjunctly irreducible and not one-generated. For every $e\in E$, there would be a proper subgraph $\G_e\leq\G$ with $e\in\G_e$. Thus, the family $(\G_e)_{e\in E}$ yields a nontrivial conjunct representation of $\G$.\\
"$\Leftarrow$": Let $\G=\evr$ be one-generated. For all $\G\sub=(E\sub,V\sub,g\sub)\leq\G$, we have that $e\in E\sub$ implies $E\sub=E$. Hence, for every conjunct representation $(\phi\hi:\G\hi\rightarrow\G)\ui$, there must be an $i\in I$ with $e \in \phi\hi[\G\hi]$. Since $e\in\phi\hi[\G\hi]\subseteq\G$, we conclude that $\phi\hi[\G\hi]=\G$ holds.
\end{proof}
\end{prop}
Dual to Birkhoff's subdirect representation theorem for algebras, we get:
\begin{satz} If $\F$ preserves arbitrary intersections, then every $\F$-graph $\G$ is a conjunct sum of conjunctly irreducible $\F$-graphs. 
\begin{proof}
The type functor $\F$ preserves arbitrary intersections. Hence, edge induced subgraphs exist (see \ref{intersection}). Obviously, $(\inkl{\G\ee}{\G}:\G\ee\hookrightarrow\G)_{e\in E}$ yields a conjunct representation of $\G$ and each $\G\ee$ is one-generated.
\end{proof}
\end{satz}

\begin{numbem}
It follows from \cite[Proposition 4.5.2]{borceux2008handbook} that the conjunctly irreducible $\F$-graphs define a family of generators. 
\end{numbem}

\begin{numbsp}
For $\F V=\mathfrak{P}_{_{1,2}} V$, the graphs $l$\footnote{The graph $l$ consist out of one vertex with a single loop.} and $K_2$\footnote{One edge with two distinct nodes.} are conjunctly irreducible. 
\end{numbsp} 

\subsection{Graph Transformations}

We want to transform graphs of type $\F_1$ into graphs of type $\F_2$. A \emph{graph transformation} is a functor $\Tt:\FGROne\rightarrow\FGRTwo$. First, we consider graph transformations induced by a natural transformation $\tau:\F_1\Rightarrow\F_2$.

\begin{center}
\begin{tikzpicture}[shorten >=0,shorten <=0,scale=1.0,state/.style={rectangle,rounded corners,draw=white,minimum height=2em,inner sep=2pt,text centered},->,>=stealth',descr/.style={fill=white,inner sep=2.5pt}]  
\node[state] (a){$E$};
\node[state, right of=a,node distance=2.12cm,anchor=center](b){$\F_1 V$};
\node[state, right of=b,node distance=3.14cm,anchor=center](c){$E$};
\node[state, right of=c,node distance=2.12cm,anchor=center](d){$\F_1 V$};
\node[state, right of=d,node distance=2.12cm,anchor=center](e){$\F_2 V$};
\path[->,font=\scriptsize](a) edge node[auto] { $g$ } (b);
\path[shorten >=16,shorten <=16,|->,font=\scriptsize](b) edge node[auto] {$$ } (c);
\path[->,font=\scriptsize](c) edge node[auto] {$g$} (d);
\path[->,font=\scriptsize](d) edge node[auto] {$\tau_V$} (e);
\end{tikzpicture}
\end{center}

The graph $\G=\evr$ is mapped to $\G=(E,V,\tau_V\circ g)$ and the naturality of $\tau$ assures that morphisms are mapped to morphisms. If $\tau$ is surjective, then $\Tt$ is surjective on objects:\\
For that, let $\tau\vv\hme:\F_2V\rightarrow\F_1 V$ be a right-inverse of $\tau_V$. Given an $\F_2$-graph $\G=\evr$, we define $\tilde{\G}:=(E,V,\tau\vv\hme\circ g)$. It holds that $\Tt(\tilde{\G})=\G.$ (see \cite[Lemma 2.3]{GummBounded})

\begin{numbsp}\label{uncoloring}

We consider colored graphs of type $X\ee\times\F(V\times X\vv)$. A natural transformation $\tau:X\ee\times\F(V\times X\vv)\Rightarrow\F V$ is given by $\tau_V=\F(\pi_V)\circ\pi_{\F(V\times X\vv)}$. This leads to a graph transformation from colored to uncolored graphs. A special case are directed graphs, because we can interpret a directed graph as an edge colored graph of type $V\times\F V$. Each edge gets its target node assigned to. The above graph transformation maps a directed graph to the underlying undirected graph.\\ 
If we model directed graphs via $\F V=V\times V$, we define $\tau_V:(v_1,v_2)\mapsto\{v_1,v_2\}$ and get a similar graph transformation from directed to undirected graphs. 
\end{numbsp}
\begin{numbsp}
If $\F$ is $M$-small\footnote{Given a set $M$, a functor $\F:\Set\rightarrow\Set$ is called $M$-\emph{small} if 
$\F X=\bigcup\{\F(\phi)[\F M]\mid\phi: M\rightarrow X\}$ for every set $X\neq\emptyset$.}, there exists a set $C$ and a surjective natural transformation $\tau$ from $C\times(-)^M$ to $\F$ (see \cite{GummBounded}). For instance, every graph of type $\mathfrak{P}_{_{1,2}}$ arises from a graph of type $(-)^{\{1,2\}}$. In this model, each edge is represented by a map, which selects the respective nodes from $V$.\\ Therefore, a simple directed graph can be considered as a family of maps $g_i:\{1,2\}\rightarrow V$ for $i\in I$ , together with the vertex set $V$. A graph $\G=((g_i)\ui,V)$ can be extended with a coloring $\gamma:V\rightarrow X\vv$ and we get $\gamma\circ\G:=((\gamma\circ g_i)\ui,V)$. If for all $i\in I$ the map $\gamma\circ g_i$ is injective, then $\gamma$ is a proper vertex coloring.
\end{numbsp}
A more general type of graph transformation is induced by $\tau:\Tt\ee\circ\F_1\Rightarrow\F_2\circ\Tt\vv$, where $\Tt\ee$ and $\Tt\vv$ are $\Set$ endofunctors. These functors manipulate the edge- and vertex set. 

\begin{center}
\begin{tikzpicture}[shorten >=0,shorten <=0,scale=1.0,state/.style={rectangle,rounded corners,draw=white,minimum height=2em,inner sep=2pt,text centered},->,>=stealth',descr/.style={fill=white,inner sep=2.5pt}]  
\node[state] (a){$E$};
\node[state, right of=a,node distance=2.12cm,anchor=center](b){$\F V$};
\node[state, right of=b,node distance=3.14cm,anchor=center](c){$\Tt\ee E$};
\node[state, right of=c,node distance=2.12cm,anchor=center](d){$\Tt\ee\F_1 V$};
\node[state, right of=d,node distance=2.12cm,anchor=center](e){$\F_2\Tt\vv V$};
\path[->,font=\scriptsize](a) edge node[auto] { $g$ } (b);
\path[shorten >=16,shorten <=16,|->,font=\scriptsize](b) edge node[auto] {$$ } (c);
\path[->,font=\scriptsize](c) edge node[auto] {$\Tt\ee(g)$} (d);
\path[->,font=\scriptsize](d) edge node[auto] {$\tau_V$} (e);
\end{tikzpicture}
\end{center}

The graph $\G=\evr$ is mapped to $\G=(\Tt\ee E,\Tt\vv V,\tau_V\circ \Tt\ee(g))$. 
\begin{numbsp}
We are not aware of an example which fully exploits this general type of transformation. We define $\Tt\ee(-):=(-)\times(-)$ and interpret it as a doubling of edges. For graphs of type $\mathfrak{P}_{_{1,2}}$, a natural transformation $\tau:\mathfrak{P}_{_{1,2}}\times\mathfrak{P}_{_{1,2}}\Rightarrow\mathfrak{P}_{_{1,2}}\circ\,\mathfrak{P}_{_{1,2}}\Rightarrow\mathfrak{P}_{_{1,2}}$ can be defined. Hence, $e\mapsto\{v_1,v_2\}$ is transformed to $(e,e)\mapsto(\{v_1,v_2\},\{v_1,v_2\})\mapsto\{\{v_1,v_2\}\}\mapsto\{v_1,v_2\}$.  
\end{numbsp}

A graph transformation which does not arise from a natural transformation is the \emph{simplification} of $\F$-graphs.
The edge set $E$ is mapped to $g[E]$, the image under $g$, and an injective structure map for the simplified graph is given by the inclusion map $\inkl{g[E]}{\F V}$. 
\begin{figure}[H]
\begin{itemize}
\item
An $\F$-graph $G$ is mapped to the simplified graph of type $F$.\\
 \begin{center}
\begin{minipage}[h]{4 cm}
\begin{tikzpicture}[>=stealth',shorten >=1pt,shorten <=1pt,state/.style={rectangle,rounded corners,draw=white,minimum height=2em,inner sep=2pt,text centered}] 	      	
 	  \node[state] (a){$E$};
 	  \node[state,below of=a,node distance=1.7cm,anchor=center](b){$\F V$};
 	  \node[state,right of=a,node distance=3cm,anchor=center](c){$ g[E]$};
 	  \node[state,below of=c,node distance=1.7cm,anchor=center](d){$\F V$};
 	  \node[below of=a,node distance=0.75cm,anchor=center](e){};
 	  \node[below of=c,node distance=0.75cm,anchor=center](f){};  	    
 	  \path[->,font=\scriptsize](a) edge  node[auto,swap] {$g$} (b); 
 	  \path[left hook->,font=\scriptsize](c) edge  node[auto] {$\inkl{g[E]}{\F V}$} (d);
 	  \path[|->,font=\scriptsize,shorten >=18pt,shorten <=18pt](e) edge  node[auto,swap] {$$} (f); 	  
\end{tikzpicture}
\end{minipage}
\end{center}
\end{itemize}
\end{figure}
\begin{figure}[H]
\begin{itemize}
\item A homomorphism
 $\phi:\G\he\rightarrow\G\hz$ is transformed into a morphism between the simplified graphs. The new edge map is the unique diagonal in the E-M-square induced by $g'$, $\F(\phi\vv)\circ\inkl{g\he[E\he]}{\F V\he}$, $g'\circ\phi\ee$ and $\inkl{g\hz[E\hz]}{\F V\hz}$.

\begin{minipage}[h]{4cm}
\begin{center}
\begin{tikzpicture}[>=stealth',shorten >=1pt,shorten <=1pt,state/.style={rectangle,rounded corners,draw=white,minimum height=2em,inner sep=2pt,text centered}] 	      	
 	  \node[state] (a){$E\he$};
 	  \node[state,below of=a,node distance=1.7cm,anchor=center](b){$\F V\he$};
 	  \node[state,right of=a,node distance=2.72cm,anchor=center](c){$E\hz$};
 	  \node[state,below of=c,node distance=1.7cm,anchor=center](d){$\F V\hz$};  	    
 	  \path[->,font=\scriptsize](a) edge  node[auto,swap] {$g\he$} (b); 
 	  \path[->,font=\scriptsize](c) edge  node[auto] {$g\hz$} (d);
 	  \path[->,font=\scriptsize](a) edge  node[auto] {$\phi\ee$} (c);
 	  \path[->,font=\scriptsize](b) edge  node[auto,swap] {$\F(\phi\vv)$} (d); 	  
\end{tikzpicture}
\end{center}
\end{minipage}
\begin{minipage}[h]{2cm}
\begin{tikzpicture}[>=stealth',shorten >=2pt,shorten <=2pt,state/.style={rectangle,rounded corners,draw=white,minimum height=2em,inner sep=2pt,text centered}]
\node[state] (a){$$};
\node[state,right of=a,node distance=1.3cm,anchor=center](b){$$};
\node[state,left of=a,node distance=0.6cm,anchor=center](c){$$};
\path[|->,font=\scriptsize](a) edge  node[auto,swap] {$$}(b);
\end{tikzpicture}
\end{minipage}
\begin{minipage}[h]{4cm}
\begin{tikzpicture}[>=stealth',shorten >=2pt,shorten <=2pt,state/.style={rectangle,rounded corners,draw=white,minimum height=2em,inner sep=2pt,text centered}]
 	  \node[state] (a){$g\he[E\he]$};
 	  \node[state,below of=a,node distance=1.7cm,anchor=center](b){$\F V\he$};
 	  \node[state,right of=a,node distance=3.72cm,anchor=center](c){$g\hz[E\hz]$};
 	  \node[state,below of=c,node distance=1.7cm,anchor=center](d){$\F V\hz$};  	    
 	  \path[left hook->,font=\scriptsize](a) edge  node[auto] {$\inkl{g\he[E\he]}{\F V\he}$} (b); 
 	  \path[left hook->,font=\scriptsize](c) edge  node[auto] {$\inkl{g\hz[E\hz]}{\F V\hz}$} (d);
 	  \path[dashed,->,font=\scriptsize](a) edge  node[auto] {$$} (c);
 	  \path[->,font=\scriptsize](b) edge  node[auto,swap] {$\F(\phi\vv)$} (d);
 	  
\end{tikzpicture}
\end{minipage}
\end{itemize}
\end{figure}

\begin{numbem}\label{coloredgraphexample}
Every coloring $\gamma:\Uu\evr\rightarrow X$ uniquely defines the colored graph $(E,V,\tilde{g})$, where $\tilde{g}:=\gamma\ee\times(\F(\id_V\times\gamma\vv)\circ g)$ is the unique mediating morphism in the diagram below. Because homomorphisms are not preserved, this construction is not functorial and does not induce a graph transformation. 
\begin{center}
\begin{tikzpicture}[state/.style={rectangle,rounded corners,draw=white,minimum height=2em,
           inner sep=2pt,text centered},->,>=stealth',descr/.style={fill=white,inner sep=2.5pt},transform shape,scale=0.8] 
					
\node[state] (a){$E$};	
\node[state, right of=a,node distance=3cm,anchor=center](b){$X\ee\times\F(V\times X\vv)$};
\node[state, right of=b,node distance=2cm,anchor=center](c){};
\node[state, below of=c,node distance=1.5cm,anchor=center](d){$\F(V\times X\vv)$};
\node[state, above of=c,node distance=1.5cm,anchor=center](e){$X\ee$};	
\node[state, right of=d,node distance=5cm,anchor=center](f){};	
\node[state, below of=f,node distance=1.5cm,anchor=center](g){$\F V$};
\node[state, above of=f,node distance=1.5cm,anchor=center](h){$\F\ X\vv$};	
\path[dashed,->,font=\scriptsize](a) edge node[descr] { $\tilde{g}$ } (b);
\path[->,font=\scriptsize](a) edge[bend left=15] node[auto] { $\gamma\ee$ } (e);
\path[->,font=\scriptsize](a) edge[bend right=20] node[auto,swap] { $g$ } (g);
\path[->,font=\scriptsize](g) edge node[auto,swap] { $\F(\gamma\vv)$ } (h);
\path[->,font=\scriptsize](b) edge node[auto] {} (e);
\path[->,font=\scriptsize](b) edge node[auto] {} (d);
\path[->,font=\scriptsize](d) edge[bend right=10] node[auto] {} (g);
\path[dashed,->,font=\scriptsize](g) edge[bend right=10] node[auto,swap] {$\F(\id_V\times\gamma\vv)$} (d);
\path[->,font=\scriptsize](d) edge node[auto] {} (h);
\path[->,font=\scriptsize](g) edge[loop right] node[auto]{$\F(\id_V)$} (g);
\end{tikzpicture}
\end{center}
\end{numbem}

From remark \ref{coloredgraphexample}, the following question arose:\\ Let $\G\he$ and $\G\hz$ be two $\mathfrak{P}$-graphs together with a homomorphism $\phi:\G\he\rightarrow\G\hz$. If an orientation\footnote{For a graph of type $\mathfrak{P}$, an \emph{orientation} is a map $\omega:E\rightarrow V$, such that $\omega(e)\in g(e)$. Also note that $(E,V,\omega)$ defines a graph of type $id$.} $\omega\he: E\he\rightarrow V\he$ is given, can we find an orientation $\omega\hz:E\hz\rightarrow V\hz$, such that $\phi$ is a homomorphism between the directed graphs $(E\he,V\he,\omega\he\times g\he)$ and $(E\hz,V\hz,\omega\hz\times g\hz)$? Generally, this is not possible, but the following holds: 

\begin{satz}
As above, let $\phi:\G\he\rightarrow\G\hz$ and an orientation $\omega\hz:E\hz\rightarrow V\hz$ be given. There exists an orientation $\omega\he: E\he\rightarrow V\he$, such that $\phi:(E\he,V\he,\omega\he\times g\he)\rightarrow(E\hz,V\hz,\omega\hz\times g\hz)$ is a homomorphism between the directed graphs.
\begin{proof}
We choose  $e\in E\he$ and define $W:=g\he(e)$. It holds that $\mathfrak{P}(\phi\vv)(g\he(e))=\phi\vv{[}g\he(e){]}$.
Because $\omega\hz$ is an orientation, it follows that $\omega\hz(\phi\ee(e))\in g\hz(\phi\ee(e))$. Furthermore, $\phi$ is a homomorphism and thus:
\[g\hz(\phi\ee(e))=\mathfrak{P}(\phi\vv)(g\he(e))=\phi\vv{[}g\he(e){]}.\]
Hence, $\omega\hz(\phi\ee(e))\in\phi\vv{[}g\he(e){]}$. Therefore, some $v\in V\he$, such that $\omega\hz(\phi\ee(e))=\phi\vv(v)$, exists. We define $\omega\he(e):=v\in g\he(e)$. Thus, $\phi$ is a homomorphism between the directed graphs. 
\end{proof}
\end{satz}

The last transformation that we mention is the minimization of $\F$-graphs. An $\F$-graph $\G$ is minimal if the only congruence on $\G$ is the diagonal relation $(\triangle_{\scriptscriptstyle{E}},\triangle_{\scriptscriptstyle{V}})$.\\ 
For any graph $\G$, let $\nabla(\G)$ be the factor graph obtained by factoring with the largest congruence relation. This congruence is the kernel of the unique homomorphism into the terminal $\F$-graph. As in \cite{Gumm_onminimal}, it can be shown that $\nabla$ is a functor from $\FGR$ to the subcategory of minimal $\F$-graphs and that minimal $\F$-graphs define an epi-reflective subcategory.\\

\section{Co-Varieties And Co-Birkhoff Theorems}\label{CoBirkSection}
In \cite{HughesCoBikhoff}, Co-Birkhoff like theorems are analyzed in an abstract categorical manner. We will use these results as a basis, in order to develop Co-Birkhoff like theorems for $\FGR$. 
\subsection{Co-Varieties}

\begin{defi} Let $\mathcal{K}$ be a class of $\F$-graphs. We define the following classes: 
\begin{itemize}
\item $\Hom(\mathcal{K})$: the class of all homomorphic images of objects from $\mathcal{K}$, 
\item $\HomM(\mathcal{K})$: the class of all homomorphic preimages of objects from $\mathcal{K}$,
\item $\Sim(\mathcal{K})$: the class of all subgraphs of objects from $\mathcal{K}$,
\item $\Sigma(\mathcal{K})$: the class of all coproducts of objects from $\mathcal{K}$.
\end{itemize}
A	class $\mathcal{K}$ is \emph {closed under $\Hom$, $\HomM$, $\Sim$ or $\Sigma$}, provided that $\Hom(\mathcal{K})\subseteq\mathcal{K}$, $\HomM(\mathcal{K})\subseteq\mathcal{K}$, $\Sim(\mathcal{K})\subseteq\mathcal{K}$, or $\Sigma(\mathcal{K})\subseteq\mathcal{K}$.
\end{defi} 

\begin{defi} A \emph{co-variety} is a class $\Kk$ of $\F$-graphs, which is closed under $\Sim$, $\Hom$ and $\Sigma$. A \emph{co-quasivariety} is a class, closed under $\Hom$ and $\Sigma$. And a \emph{complete co-variety} is class, closed under $\HomM$, $\Sim$, $\Hom$ and $\Sigma$. 
\end{defi}

\begin{numbem} Let $\Kk$ be a class of $\F$-graphs. It can be shown that $\Sim\Hom\Sigma(\Kk)$ is the smallest co-variety containing $\Kk$.
\end{numbem}

From the previous sections, we know that $\FGR$ has coproducts (theorem \ref{ColimitTheorem}), has an epi-regular mono-factorisation system (theorem \ref{FGRfactorization}), is regular-well-powered (theorem \ref{RegularmonoCharacterization}), has enough regular injectives (lemma \ref{ExtensionCoFree},\ref{unitembedding}), has binary products (theorem \ref{products}) and a terminal object (theorem \ref{cofreegraphcharacterization}). Thus, we can apply \cite[Theorem 2.3]{HughesCoBikhoff} and \cite[Theorem 3.6]{HughesCoBikhoff}. Furthermore, due to theorem \ref{RegInjChar}, it is sufficient to consider cofree $\F$-graphs instead of arbitrary regular-injective $\F$-graphs. 

\begin{defi} 
Let $\G=\evr$ be an $\F$-graph and $\Cc X$ the cofree graph over a color set $X=(X\ee,X\vv)$. We define $\col_X(\G)=\{\gamma\mid\gamma:\Uu\G\rightarrow X\}$ to be the collection of all colorings of $\G$. A subset $P\subseteq(X\ee,X\vv)$ is called \emph{pattern} over $\Cc X$. We say that a pattern $P$ holds in $\G$ if for all $\gamma\in\col_X(\G)$, we have that $\overline{\gamma}[\G]\leq\hat{P}:=[P]\ug$ (recall definition \ref{cogensubgrdefi}). For that, we  write $\G\vDash P$. This means that for all colorings, the induced homomorphism $\overline{\gamma}$ factors through the largest in $P$ contained subgraph, or in other words that $\G$ is regular-projective with respect to the inclusion morphism $\inkl{\hat{P}}{\C X}$.
\begin{figure}[H]
\begin{center}
\begin{tikzpicture}[description/.style={fill=white,inner sep=2pt},>=stealth']
\matrix (m) [matrix of math nodes, row sep=3.5em,
column sep=2.5em, text height=1.3ex, text depth=0.25ex] 
{\G & & \Cc X\\
    & & \hat{P} \\ };
\path[->,font=\scriptsize](m-1-1) edge node[auto] {$ \overline{\gamma} $} (m-1-3);
\path[dashed,->,font=\scriptsize](m-1-1) edge node[auto,swap] {$ $} (m-2-3);
\path[right hook->,shorten <=4,font=\scriptsize](m-2-3) edge node[auto,swap] {$\inkl{\hat{P}}{\Cc X}$} (m-1-3);
\end{tikzpicture}
\end{center}
\end{figure}
Let $\Pp$ be a collection of patterns and $\Kk$ a class of $\F$-graphs. We say $\G\vDash \Pp$ if every $P\in \Pp$ holds in $\G$. Analogously, $\Kk\vDash P$ provided $\G\vDash P$ for every $\G\in\Kk$.
\end{defi}

\begin{satz}\footnote{\cite[Theorem 3.6.21]{HughesDiss}}\label{covarietychar} A class $\Kk$ of $\F$-graphs is a co-variety if and only if there is a collection $\Pp$ of patterns, such that for all $\G$, it holds that $\G\in\Kk\Leftrightarrow\G\vDash \Pp.$
\end{satz}

We define: $\Graph(\Pp):=\{\G\in\FGR\mid \G\vDash\Pp\}$ and $\Pat(\Kk):=\{P \mid \Kk\vDash P\}.$

\begin{satz}[Co-Birkhoff Theorem For $\F$-Graphs]\footnote{\cite[Theorem 2.3]{HughesCoBikhoff}}\label{CoBiGeneral} For a class $\Kk$ of $\F$-graphs, we have: $\Sim\Hom\Sigma(\Kk)=\Graph(\Pat(\Kk)).$
\end{satz}

With an additional condition, we can assure that every co-variety can be specified by a set of patterns. 
\begin{defi}
An $\F$-graph $\G=\evr$ is bounded by $X\vv$ if for every $e\in E$ the edge induced subgraph $\G\ee=(E\ee,V\ee,g\ee)$ has at most $|X\vv|$ nodes, \textit{i.e}, $|V\ee|\leq|X\vv|$. 
\end{defi}

For instance, graphs of type $\mathfrak{P}_{_{1,2}}$ are bounded by $\{1,2\}$.
For fixed $X$, we define $\Pat_X(\Kk):=\{P\subseteq X\mid \Kk\vDash P\}.$

\begin{satz}[Co-Birkhoff Theorem For Bounded $\F$-Graphs]\label{CoBirkTh} Let every $\G\in\FGR$ be bounded by $X\vv$ and let $\F$ preserve arbitrary intersections. A co-variety can be specified by a set of patterns over some cofree $\F$-graph $\Cc X$. That is, for an arbitrary class of $\F$-graphs $\Kk$, we have:\[\Sim\Hom\Sigma(\Kk)=\Graph(\Pat_X(\Kk)).\]
\begin{proof}
From \ref{CoBiGeneral}, the "$\subseteq$"-direction follows. For the converse, let $\G\in\Graph(\Pat_X(\Kk))$ and $X:= (\{\ast\},X\vv)$. Because $\G$ is bounded, we have that for every edge induced subgraph $\G\ee\leq\G$, there is an injective homomorphism $\G\ee\hookrightarrow\Cc(\{\ast\},X\vv)$. Hence, $\G$ is the conjunct sum of subgraphs from $\Cc X$ and there is a surjective homomorphism $\sum\G\ee\twoheadrightarrow\G$. We show that every $\G\ee$ is contained in $\tilde{\Kk}:=\Sim\Hom\Sigma(\Kk)$.\\ 
Let $\hat{P}:=\bigcup \{\overline{\gamma}[\G]\mid \G\in\tilde{\Kk},\;\gamma\in\col_{(\{\ast\},X\vv)}(\G)\}$. Because $\hat{P}$ is a conjunct sum of homorphic images from elements of $\tilde{\Kk}$, it holds that $\hat{P}\in\tilde{\Kk}$. Furthermore, we know that $\G\in\Graph(\Uu \hat{P})$ and due to the regular-injectivity of $\Cc X$, every $\G\ee$ is contained in $\Graph(\Uu \hat{P})$. As $\G\ee$ is a subgraph of $\Cc X$, it follows that $\G\ee\leq \hat{P}$. Thus, it follows that $\G\ee\in\tilde{\Kk}$.
\end{proof}
\end{satz}

Next, we characterize the patterns which define a co-variety. 
\begin{defi}  We say that $\G\sub\leq\G$ is \emph{invariant in} $\G$ if for every endomorphism $\phi:\G\rightarrow\G$, we have that $\phi[\G\sub]\leq\G\sub$. 
\end{defi}
The following lemma is inspired by \cite{GummCompCov}.
\begin{lemma}\label{invariancelemma}
Let $\Cc X$ be the cofree $\F$-graph over $X$ and $P\subseteq X$ a pattern in $X$. It holds that: $\hat{P}\in\Graph(P)\Longleftrightarrow \hat{P}\;\text{is invariant in}\;\Cc X.$
\begin{proof}
"$\Rightarrow$": For $\phi:\Cc X\rightarrow\Cc X$, we consider the restriction $\phi\circ\inkl{\hat{P}}{\Cc X}$ to $\hat{P}$. Because $\hat{P}\in\Graph(P)$, the restriction must factor through $\hat{P}$ and this means $\phi[\hat{P}]\leq \hat{P}$.\\
"$\Leftarrow$": Let $\phi:\hat{P}\rightarrow \Cc X$ be a homomorphism.  We can extend $\phi$ to $\psi:\Cc X\rightarrow\Cc X$ (see lemma \ref{ExtensionCoFree}). The invariance of $\hat{P}$ in $\Cc X$ yields $\phi[\hat{P}]=(\psi\circ\inkl{\hat{P}}{\Cc X})[\hat{P}]\leq\hat{P}$.
\end{proof}
\end{lemma} 

\begin{satz}\label{CovarChar} For $X\vv$-bounded $\F$-graphs, each co-variety corresponds exactly to the invariant subgraphs of $\Cc(\{\ast\},X\vv)$.
\begin{proof}
For a given co-variety $\Kk$, we define $\tilde{\G}:=\bigcup \{\overline{\gamma}[\G]\mid\G\in\Kk,\;\gamma\in\col_{(\{\ast\},X\vv)}(\G)\}$. As in the proof of theorem \ref{CoBirkTh}, we conclude that $\tilde{\G}\in\Kk$. Because of $\Kk\subseteq\Graph(\tilde{\G})$ and lemma \ref{invariancelemma}, it follows that $\tilde{\G}$ is an invariant subgraph of $\Cc(\{\ast\},X\vv)$.\\ On the other hand, if $\G\sub$ is a subgraph of $\Cc(\{\ast\},X\vv)$, it follows from theorem \ref{covarietychar} that $\Graph(\G\sub)$ is a co-variety. 
\end{proof}
\end{satz}

\begin{numbsp}
\begin{itemize}
\item[(i)]Considering graphs of type $\mathfrak{P}_{_{1,2}}$, for any color set $X=(X\ee,X\vv)$ with $|X\vv|\geq 2$, we put $\Kk=\{l_1\}$ and the resulting co-variety contains all graphs consisting exclusively out of loops. Other choices for $\Kk$ lead to the co-variety of all graphs. That is because the only pattern that holds in a graph $\G\neq l_1$ is $X$.
We also recognize that all loops or $\Cc X$ itself are the only invariant subgraphs in $\Cc X$.
\item[(ii)] We can consider graphs of type $\mathfrak{P}_{_{1,2}}$ as subgraphs of type $\mathfrak{P}$. For any graph $\G$ which is not $l_1$, the construction $\Sim\Hom\Sigma(\G)$ leads to the co-variety of all graphs of type $\mathfrak{P}_{_{1,2}}$ within the category of $\mathfrak{P}$-graphs.

\item[(iii)] We consider patterns over the terminal graph $\Cc(X\ee,X\vv)$ (see remark \ref{cofreeterminal}) in the category of $(X\ee,X\vv)$ colored graphs. For appropriate subsets $Y\ee\subseteq X\ee$, $Y\vv\subseteq X\vv$, we get the co-variety of $(Y\ee,Y\vv)$ colored graphs.\\ Also note that in the terminal graph $\Cc(X\ee,X\vv)$ all subgraphs are endomorphism invariant. 
\end{itemize} 
\end{numbsp}

All the examples above arise from subfunctors, \textit{i.e.}, $\mathfrak{P}_{_{1}}\subseteq\mathfrak{P}_{_{1,2}}$, $\mathfrak{P}_{_{1,2}}\subseteq\mathfrak{P}$ and $Y\ee\times\F((-)\times Y\vv)\subseteq X\ee\times\F((-)\times X\vv)$. As in \cite{TobiSchDiss}, it can be shown that a subfunctor always induces a co-variety. 

\subsection{Co-Quasivarieties} 
We restate conditional co-equations (\cite[Section 3]{HughesModal}) in terms of $\F$-graphs.  

\begin{defi}
A \emph{conditional pattern} over an arbitrary  $\F$-graph $\G$ is a subset $P\subseteq\Uu\G$. We say that $\tilde{\G}\vDash\ug P$ if for every homomorphism $\phi:\tilde{\G}\rightarrow\G$, it holds that $\phi[\tilde{\G}]\leq \hat{P}$.
\end{defi}

\begin{defi}
For two patterns $P,Q\subseteq X$ over $\Cc X$, we say that $P\Rightarrow Q$ holds in $\G$ if $\G\vDash P$ implies $\G\vDash Q$  and we write $\G\vDash P\Rightarrow Q$.
\begin{figure}[H]
\begin{center}
\begin{minipage}[h]{4cm}
\begin{center}
\begin{tikzpicture}[description/.style={fill=white,inner sep=2pt},>=stealth']
\matrix (m) [matrix of math nodes, row sep=3.5em,
column sep=2.5em, text height=1.3ex, text depth=0.25ex] 
{\G & & \Cc X\\
    & & \hat{P} \\ };

\path[->,font=\scriptsize](m-1-1) edge node[auto] {$ \overline{\gamma} $} (m-1-3);
\path[dashed,->,font=\scriptsize](m-1-1) edge node[auto,swap] {$ $} (m-2-3);
\path[right hook->,shorten <=4,font=\scriptsize](m-2-3) edge node[auto,swap] {$\inkl{\hat{P}}{\Cc X}$} (m-1-3);
\end{tikzpicture}
\end{center}
\end{minipage}
\hspace{0.5cm}
\begin{minipage}[h]{1cm}
$\Longrightarrow$
\end{minipage}
\hspace{0.5cm}
\begin{minipage}[h]{4cm}
\begin{center}
\begin{tikzpicture}[description/.style={fill=white,inner sep=2pt},>=stealth']
\matrix (m) [matrix of math nodes, row sep=3.5em,
column sep=2.5em, text height=1.3ex, text depth=0.25ex] 
{\G & & \Cc X\\
    & & \hat{Q}\\ };

\path[->,font=\scriptsize](m-1-1) edge node[auto] {$ \overline{\gamma} $} (m-1-3);
\path[dashed,->,font=\scriptsize](m-1-1) edge node[auto,swap] {$ $} (m-2-3);
\path[right hook->,shorten <=4,font=\scriptsize](m-2-3) edge node[auto,swap] {$\inkl{\hat{Q}}{\Cc X}$} (m-1-3);
\end{tikzpicture}
\end{center}
\end{minipage}

\end{center}
\end{figure}
\end{defi}

Obviously, if $P\subseteq Q$ and $\G\vDash P$, then $\G\vDash P\Rightarrow Q$.\\

The following lemma states an equivalence between conditional patterns and implications. 

\begin{hsatz}
Let $P$ and $Q$ be patterns over $\Cc X$. There exists an $\F$-graph $\G$ and a conditional pattern $R$ over $\G$, such that for all $\tilde{G}$ it holds that:
\[\tilde{G}\vDash P\Rightarrow Q\;\; \text{if and only if}\;\; \tilde{G}\vDash\ug R.\]
\begin{proof}
We can put $\G=\hat{P}$ and $R=\Uu\hat{P}\cap Q$. On the other hand, if we consider $\Uu\G$ and $R$ as patterns over $\Cc\Uu\G$, it holds that $\tilde{G}\vDash\ug R \Leftrightarrow \tilde{G}\vDash \Uu\G\Rightarrow P.$
\end{proof}
\end{hsatz}
A detailed proof for coalgebras can be found in \cite[Theorem 3.6.21]{HughesDiss}.
Next, we formulate the Quasi Co-Birkhoff Theorem (\cite[Corollary 3.6]{HughesModal}). 

\begin{satz}
A class $\Kk$ of $\F$-graphs is a quasi co-variety if and only if there is a collection $\Pp$ of conditional patters, such that
$\G\in\Kk \;\; \text{iff}\;\;\forall P\in\Pp:\G\vDash\ug P.$
\end{satz}
For a collection of implications $\mathcal{I}$ and some class $\Kk$ of $\F$-graphs, we define:
\[\Graph(\mathcal{I}):=\{\G\in\FGR\mid \G\vDash\mathcal{I}\}\;\text{and}\;\Imp(\Kk):=\{P\Rightarrow Q \mid  \Kk\vDash P\Rightarrow Q\}\]

\begin{folg}
For any class $\Kk$ of $\F$-graphs, we have $\Hom\Sigma(\Kk)=\Graph(\Imp(\Kk)).$
\end{folg}

\subsection{Complete Co-Varieties}
The additional closure under homomorphic preimages leads to complete co-varieties. 

\begin{satz}
From \cite[Theorem 3.6]{HughesCoBikhoff}: Let $X$ be a pattern over the terminal graph. We have that $\HomM\Sim\Hom\Sigma(\Kk)=\Graph(\Pat_X(\Kk)).$
\end{satz}

\begin{numbsp}
In the category of $(X\ee,X\vv)$ colored graphs $\Cc(X\ee,X\vv)$ is terminal. Hence, each subgraph of $\Cc(X\ee,X\vv)$ defines a complete co-variety. For instance, considering vertex colored graphs of type $\mathfrak{P}_{_{1,2}}$, the complete co-variety generated by $\Kk=\{K_n\}$ leads to all $n$-colored graphs. \\ 
\end{numbsp}
\begin{numbem}
A complete co-variety $\Kk$ is also closed under total graph relations. To see this, let $\G\he$ and $\G\hz$ be related via $\G\hr$. If $\G\he\in\Kk$, then by closure under $\HomM$ also $\G\hr\in\Kk$. The closure under $\Hom$ yields $\G\hz\in\Kk$.
\end{numbem}
\section{Acknowledgments}
For support and advice we want to thank Stefan E. Schmidt. Also special thanks go to Andreas Thom, for discussions about cofree $\F$-graphs and the construction of $\F$-graph products.

    
\begingroup
\renewcommand{\bibname}{References}
\setlength\bibitemsep{0pt}

\endgroup

\end{document}